\documentclass[12pt]{article}
\usepackage{amsmath,amssymb,amsbsy,amsfonts,amsthm,
                          latexsym,amsopn,amstext,amsxtra,euscript,amscd}
\usepackage{graphicx}
\begin{document}

\newtheorem{theorem}{Theorem}
\newtheorem{acknowledgement}[theorem]{Acknowledgement}
\newtheorem{algorithm}[theorem]{Algorithm}
\newtheorem{axiom}[theorem]{Axiom}
\newtheorem{case}[theorem]{Case}
\newtheorem{claim}[theorem]{Claim}
\newtheorem{conclusion}[theorem]{Conclusion}
\newtheorem{condition}[theorem]{Condition}
\newtheorem{conjecture}[theorem]{Conjecture}
\newtheorem{corollary}[theorem]{Corollary}
\newtheorem{criterion}[theorem]{Criterion}
\newtheorem{definition}[theorem]{Definition}
\newtheorem{example}[theorem]{Example}
\newtheorem{exercise}[theorem]{Exercise}
\newtheorem{lemma}[theorem]{Lemma}
\newtheorem{notation}[theorem]{Notation}
\newtheorem{problem}[theorem]{Problem}
\newtheorem{question}[theorem]{Question}
\newtheorem{proposition}[theorem]{Proposition}
\newtheorem{remark}[theorem]{Remark}
\newtheorem{solution}[theorem]{Solution}
\newtheorem{summary}[theorem]{Summary}

%% DEFINITIONS
\def\scr{\scriptstyle}
\def\\{\cr}
\def\({\left(}
\def\){\right)}
                      \def\[{\left[}
\def\]{\right]}
\def\<{\langle}
\def\>{\rangle}
\def\fl#1{\left\lfloor#1\right\rfloor}

\def\rf#1{\left\lceil#1\right\rceil}

\def\defn{\noindent{\bf Definition\/}\ \ }

\def\area#1{ \left|#1\right|}

\def\Z{\mathbb Z}
\def\R{\mathbb R}
\def\N{\mathbb N}
\def\cS{\mathcal S}
\def\e{{\bf e}}
\def\cR{\mathcal R}
\def\cP{\mathcal P}
\def\cQ{\mathcal Q}
\def\cT{\mathcal T}
\def\ord{{\mathrm ord\/}}
\def\eps{\varepsilon}
\def\ep{{\mathbf{\,e}}_p}
\def\epp{{\mathbf{\,e}}_{p-1}}

\def\xxx{\vskip5pt\hrule\vskip5pt}
\def\yyy{\vskip5pt\hrule\vskip2pt\hrule\vskip5pt}

\newcommand{\comm}[1]{\marginpar{\vskip-\baselineskip
%raise the marginpar a bit
\raggedright\footnotesize
\itshape\hrule\smallskip#1\par\smallskip\hrule}}

                     \newcommand{\1}{1\!{\mathrm l}} \newcommand{\croix}{
set\frown}

%%%%%%%%%%%%%%%%%%%%%%%%%%%%%%%%%%%%%%%%%%%%%%%%%%%%%%%%
%%%%%%%%%%%%%%%%%%%%%%%%%%%%%%%%%%%%%%%%%%%%%%%%%%%%%%%%
%%%%%%%%%%%%%%%%%%%%%%%%%%%%%%%%%%%%%%%%%%%%%%%%%%%%%%%%
%%%%%%%%%%%%%%%%%%%%%%%%%%%%%%%%%%%%%%%%%%%%%%%%%%%%%%%%

%%%%%%%   STANDARD STUFF %%%%%%%%%

%%%%%%%%%%%%%%%%%%%%%%%%%%%%%%%%%%%%%%%%%%%%%%%%%%%%%%%%
%%%%%%%%%%%%%%%%%%%%%%%%%%%%%%%%%%%%%%%%%%%%%%%%%%%%%%%%
%%%%%%%%%%%%%%%%%%%%%%%%%%%%%%%%%%%%%%%%%%%%%%%%%%%%%%%%
%%%%%%%%%%%%%%%%%%%%%%%%%%%%%%%%%%%%%%%%%%%%%%%%%%%%%%%%

%  use the AMS-Euler Fraktur fonts %%%%%%%%%%%%%%%%%%%%%%%%%%%%%%%%%%
\newfont{\teneufm}{eufm10}
\newfont{\seveneufm}{eufm7}
\newfont{\fiveeufm}{eufm5}
%%%%%%%%%%%%%%%%%%%%%%%%%%%%%%%%%
%
%  allow automatic size selection in math mode %
%%%%%%%%%%%%%%%%%%%%%%%%%%%%%%%%%
\newfam\eufmfam
                                         \textfont\eufmfam=\teneufm
\scriptfont\eufmfam=\seveneufm
                                         \scriptscriptfont\eufmfam=\fiveeufm
%%%%%%%%%%%%%%%%%%%%%%%%%%%%%%%%%

%
%  \frak works on a single symbol at a time...
%
\def\frak#1{{\fam\eufmfam\relax#1}}
%

%%%%%%%%%%%%%%%%%%%  bbb-matter

%
\def\bbbr{{\rm I\!R}} %reelle Zahlen
\def\bbbm{{\rm I\!M}}
\def\bbbn{{\rm I\!N}} %natuerliche Zahlen
\def\bbbf{{\rm I\!F}}
\def\bbbh{{\rm I\!H}}
\def\bbbk{{\rm I\!K}}
\def\bbbp{{\rm I\!P}}
\def\bbbone{{\mathchoice {\rm 1\mskip-4mu l} {\rm 1\mskip-4mu l} {\rm
1\mskip-4.5mu l} {\rm 1\mskip-5mu l}}}
\def\bbbc{{\mathchoice {\setbox0=\hbox{$\displaystyle\rm C$}
\hbox{\hbox to0pt{\kern0.4\wd0\vrule height0.9\ht0\hss}\box0}}
{\setbox0=\hbox{$\textstyle\rm C$}\hbox{\hbox to0pt{\kern0.4\wd0\vrule
height0.9\ht0\hss}\box0}} {\setbox0=\hbox{$\scriptstyle\rm
C$}\hbox{\hbox to0pt{\kern0.4\wd0\vrule height0.9\ht0\hss}\box0}}
{\setbox0=\hbox{$\scriptscriptstyle\rm C$}\hbox{\hbox
to0pt{\kern0.4\wd0\vrule height0.9\ht0\hss}\box0}}}}
\def\bbbq{{\mathchoice {\setbox0=\hbox{$\displaystyle\rm
Q$}\hbox{\raise 0.15\ht0\hbox to0pt{\kern0.4\wd0\vrule
height0.8\ht0\hss}\box0}}
{\setbox0=\hbox{$\textstyle\rm Q$}\hbox{\raise 0.15\ht0\hbox
to0pt{\kern0.4\wd0\vrule height0.8\ht0\hss}\box0}}
{\setbox0=\hbox{$\scriptstyle\rm
Q$}\hbox{\raise 0.15\ht0\hbox to0pt{\kern0.4\wd0\vrule
height0.7\ht0\hss}\box0}} {\setbox0=\hbox{$\scriptscriptstyle\rm
Q$}\hbox{\raise 0.15\ht0\hbox
to0pt{\kern0.4\wd0\vrule height0.7\ht0\hss}\box0}}}}
\def\bbbt{{\mathchoice {\setbox0=\hbox{$\displaystyle\rm
T$}\hbox{\hbox to0pt{\kern0.3\wd0\vrule height0.9\ht0\hss}\box0}}
{\setbox0=\hbox{$\textstyle\rm T$}\hbox{\hbox
to0pt{\kern0.3\wd0\vrule height0.9\ht0\hss}\box0}}
{\setbox0=\hbox{$\scriptstyle\rm T$}\hbox{\hbox
to0pt{\kern0.3\wd0\vrule height0.9\ht0\hss}\box0}}
{\setbox0=\hbox{$\scriptscriptstyle\rm T$}\hbox{\hbox
to0pt{\kern0.3\wd0\vrule
height0.9\ht0\hss}\box0}}}}
\def\bbbs{{\mathchoice {\setbox0=\hbox{$\displaystyle     \rm S$}
\hbox{\raise0.5\ht0\hbox to0pt{\kern0.35\wd0\vrule
height0.45\ht0\hss}\hbox to0pt{\kern0.55\wd0\vrule
height0.5\ht0\hss}\box0}}
{\setbox0=\hbox{$\textstyle        \rm S$}\hbox{\raise0.5\ht0\hbox
to0pt{\kern0.35\wd0\vrule height0.45\ht0\hss}\hbox
to0pt{\kern0.55\wd0\vrule height0.5\ht0\hss}\box0}}
{\setbox0=\hbox{$\scriptstyle      \rm S$}\hbox{\raise0.5\ht0\hbox
to0pt{\kern0.35\wd0\vrule height0.45\ht0\hss}\raise0.05\ht0\hbox
to0pt{\kern0.5\wd0\vrule height0.45\ht0\hss}\box0}}
{\setbox0=\hbox{$\scriptscriptstyle\rm S$}\hbox{\raise0.5\ht0\hbox
to0pt{\kern0.4\wd0\vrule height0.45\ht0\hss}\raise0.05\ht0\hbox
to0pt{\kern0.55\wd0\vrule height0.45\ht0\hss}\box0}}}}
\def\bbbz{{\mathchoice {\hbox{$\sf\textstyle Z\kern-0.4em Z$}}
{\hbox{$\sf\textstyle Z\kern-0.4em Z$}} {\hbox{$\sf\scriptstyle
Z\kern-0.3em Z$}} {\hbox{$\sf\scriptscriptstyle Z\kern-0.2em Z$}}}}
\def\ts{\thinspace}

\def\squareforqed{\hbox{\rlap{$\sqcap$}$\sqcup$}}
\def\qed{\ifmmode\squareforqed\else{\unskip\nobreak\hfil
\penalty50\hskip1em\null\nobreak\hfil\squareforqed
\parfillskip=0pt\finalhyphendemerits=0\endgraf}\fi}%%

%%%%%%%%%%%%%%%%%%%%%%%%%
% Alphabet calligraphic %
%%%%%%%%%%%%%%%%%%%%%%%%%
\def\cA{{\mathcal A}}
\def\cB{{\mathcal B}}
\def\cC{{\mathcal C}}
\def\cD{{\mathcal D}}
\def\cE{{\mathcal E}}
\def\cF{{\mathcal F}}
\def\cG{{\mathcal G}}
\def\cH{{\mathcal H}}
\def\cI{{\mathcal I}}
\def\cJ{{\mathcal J}}
\def\cK{{\mathcal K}}
\def\cL{{\mathcal L}}
\def\cM{{\mathcal M}}
\def\cN{{\mathcal N}}
\def\cO{{\mathcal O}}
\def\cP{{\mathcal P}}
\def\cQ{{\mathcal Q}}
\def\cR{{\mathcal R}}
\def\cS{{\mathcal S}}
\def\cT{{\mathcal T}}
\def\cU{{\mathcal U}}
\def\cV{{\mathcal V}}
\def\cW{{\mathcal W}}
\def\cX{{\mathcal X}}
\def\cY{{\mathcal Y}}
\def\cZ{{\mathcal Z}}%%

\newcommand{\li}{\operatorname{li}}

%%%%%%%%%%%%%%%%%%%%%%%%%%%%%%%%%%%%%%%%
%%%%%%%%%%%%%%%%%%%%%%%%%%%%%%%%%%%%%%%%
%%%%%%%  END OF STANDARD STUFF %%%%%%%%%
%%%%%%%%%%%%%%%%%%%%%%%%%%%%%%%%%%%%%%%%
%%%%%%%%%%%%%%%%%%%%%%%%%%%%%%%%%%%%%%%%

\def\cNt{\widetilde{\cN}}

\title{On the Convex Closure of \\ the Graph of  Modular Inversions}

\author{
{\sc Mizan R. Khan}
\\ Department of Mathematics and Computer Science \\Eastern
Connecticut State University
\\ Willimantic, CT 06226, USA \\ {\tt khanm@easternct.edu}
\and
              {\sc Igor E. Shparlinski} \\ Department of
Computing\\ Macquarie University\\ Sydney, NSW 2109, Australia\\{\tt
igor@ics.mq.edu.au}  \and   {\sc Christian L.
Yankov}\\ Department of Mathematics and Computer Science \\Eastern
Connecticut State University \\Willimantic, CT
06226, USA\\{\tt yankovc@easternct.edu} }
\date{November 21, 2007}
\maketitle

\begin{abstract}
In this paper we give upper and lower bounds as well as a heuristic
estimate on the number of vertices of the convex
closure of the set $$ G_n=\left\{  (a,b)\ : \ a,b\in \Z, ab \equiv 1
\pmod{n},\ 1\leq a,b\leq n-1\right\}.
$$
The heuristic is based on an asymptotic formula of R\'{e}nyi and
Sulanke.
After describing two algorithms to determine the convex closure, we
compare the numeric results with the heuristic estimate. The numeric
results do not agree with the heuristic estimate --- there are some
interesting peculiarities for which we provide a
heuristic explanation. We then describe some numerical work on the 
convex closure of the
graph of random quadratic and cubic polynomials over $\Z_n$. In this
case the numeric results are in much closer agreement with the
heuristic, which strongly suggests that the the curve $xy=1\pmod{n}$
is ``atypical''.
\end{abstract}

\section{Introduction}

Let $G_n$ be the set
$$ G_n=\left\{  (a,b)\ : \ a,b\in \Z, ab \equiv 1 \pmod{n},\ 1\leq
a,b\leq n-1\right\},$$
whose cardinality is given by the Euler function $\varphi(n)$.
If we scale by a factor of $1/n$ we get the set of points
                      $n^{-1}G_n$, which is uniformly distributed in the unit
square.
More precisely, if $\Omega \subseteq [0,1]^2$ has piecewise
smooth boundary and $N(\Omega,n)$ is the cardinality of the
intersection $\Omega \cap n^{-1}G_n$, then it is natural
to expect, and in fact can be proved by using the bounds of {\it
Kloosterman sums\/}, that \begin{equation}
\label{eq:Asymp} \left |\area
\Omega-\frac{N(\Omega,n)}{\varphi(n)}\right | \to 0 \quad \text{as} \
n \to \infty,
\end{equation} where $\area \Omega$ is the area of $\Omega$.
Figure~1, generated by {\sc Maple},
illustrates this property.

\centerline{ \includegraphics*[scale=.5]{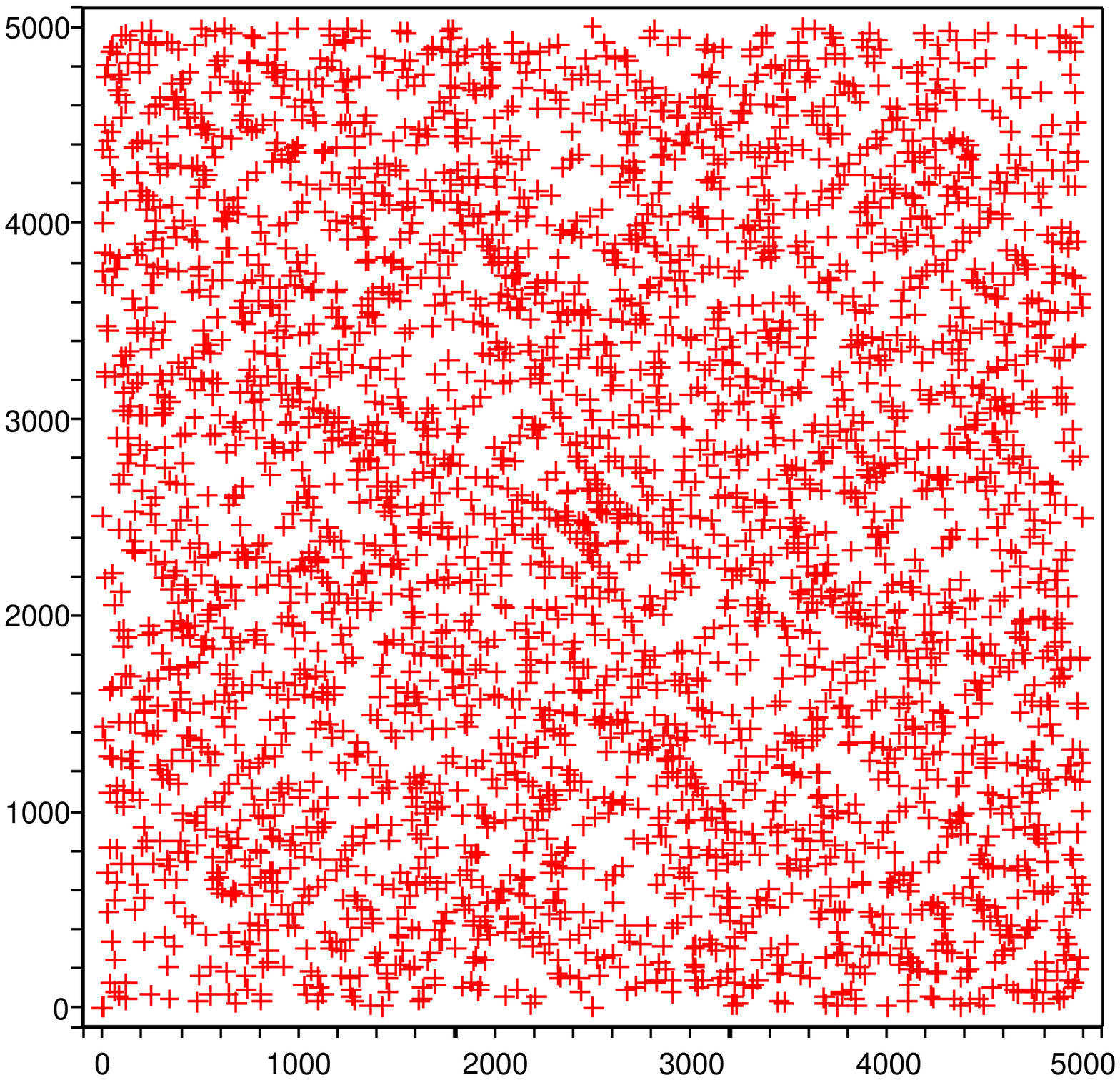} }
\centerline{{\bf Fig. 1.} The graph $G_{5001}$} \medskip
\medskip

Quantitative forms of~\eqref{eq:Asymp} have been given in a number of
works,
see~\cite{CZ,GrShZa,VZ,Zha,Zhe} and references therein.
For example, it follows from more general results of~\cite{GrShZa} that
for
primes $p$, \begin{equation} \label{eq:CZ-Bd} \left |\area
\Omega-\frac{N(\Omega,p)}{p-1}\right | =O\(p^{-1/4}\log
p\right ), \end{equation} where the implied constant depends only on
$\Omega$.

Here  we continue to study some geometric properties of the set $G_n$
and
in particular concentrate on the convex closure $C_n$ of   $G_n$.
One of our questions of interest is the behavior  of $v(n)$ and
$V(N)$, where $v(n)$ denotes the number of vertices of $C_n$ and
$V(N)$ denotes the average,
$$V(N)=\frac{1}{N-1}\sum_{n=2}^N v(n). $$

We demonstrate that the theoretic and algorithmic study of $v(n)$
has surprising links with various areas of number theory, such as
bounds of exponential sums, distribution of
divisors of ``typical'' integers and integer factorisation. On the
other hand, we present heuristic estimates $h(n)$ and $H(N)$ for
$v(n)$ and $V(N)$, respectively.  These heuristic estimates arise by
viewing
$G_n$ as a set of points that are randomly distributed and then using
the
result of R\'{e}nyi and Sulanke~\cite[Satz~1]{ReSu1}.
On comparing with our numeric results we see that although the
heuristic prediction $H(N)$ gives an adequate idea about the type of
growth of $V(N)$,
there is a deviation which behaves quite regularly and thus probably
reflects certain
other hidden effects.
We suggest some explanation. We also examine
numerically some other interesting peculiarities in the behaviour
of $v(n)$ which lead us to several open questions.

Finally, we present some numerical evidence suggesting that
the above effects do not arise for sets of points on other curves
which behave more like truly random sets of points,
which makes the study of $G_n$ even more interesting.

We note that some other geometric properties of the points of $G_n$
have recently been considered in~\cite{ShpWin}.
A survey of recent results about the distribution of points
of $G_n$ and more general sets corresponding to
congruences of the type  $ab\equiv \lambda \pmod n$
with some fixed $\lambda$, are given in~\cite{Shp}.

\section{Some Preliminary Observations}

\subsection{General structure of $C_n$}
We begin with a simple (but useful) remark on two lines of symmetry of
$G_n$.

\begin{proposition}\label{Sym-Lem}
The points of $G_n$ are symmetrically distributed about the lines
$y=x$ and $x+y=n$.
\end{proposition}

Therefore, if $\(  a,b\)  \in G_n$, then its reflection in $y=x$,
$(b,a)$,  and its reflection in $x+y=n$,
$\(n-b,n-a\)  $, are elements of $G_n$. Consequently, $\(  a,b\)  $
is a boundary point of $C_n$, if and only if $\(
b,a\)  ,(n-b,n-a)$ and $(n-a,n-b)$ are boundary points of $C_n$.

Our next result shows that  $C_n$ is always a convex polygon with
nonempty interior, except
when $n=2,3,4,6,8,12,24$.

\begin{proposition}
$\area{C_n}  =0$, if and only if $n=2,3,4,6,8,12$ or $24$.
\end{proposition}

\begin{proof}
This follows by observing that for these moduli all of the elements
in $\Z_{n}^{\ast}$
(that is, all units of the residue ring modulo $n$)
have order 2.
Consequently, for these moduli all of the elements of $G_n$ lie on
the line $y=x$.
\end{proof}

{From} now on we  typically exclude the cases $n=2,3,4,6,8,12$ and $24$.

\subsection{Points in the triangle $\cT_n$} \label{sec:Triangle Points}

By Proposition~\ref{Sym-Lem} we only need to know the vertices of $C_n$
that
lie in the triangle $\cT_n$ with vertices $(0,0),(0,n)$ and $\(
n/2,n/2\) $, to determine $C_n$. We denote the
vertices of $C_n  $ that lie in the triangle $\cT_{n}$ by $$
(a_0,b_{0}),(a_1,b_{1}),\ldots,(a_s,b_{s}) \in C_n \cap \cT_n, $$
where $a_{0}<a_{1}<\ldots<a_{s}$.

\begin{proposition}\label{ver-cond}
We have the following:
\begin{enumerate}
\item $(a_0,b_0)=(1,1)$;
\item $a_{i}<b_{i}$ for $i=1,\ldots ,s$;
\item $b_{0}<b_{1}<\ldots <b_{s}$.
\item $b_i-a_i < b_{i+1}-a_{i+1}$ for $i=0,\dots,s-1$.
\end{enumerate}
\end{proposition}
\begin{proof}
Assertions 1 and 2 are clear. Assertions 3 and 4 follow from the
following observation.
The line through $(a_i,b_i)$ and its symmetric counterpart
$(n-b_i,n-a_i)$ intersects the line $x+y=n$ at the point
$((n-b_i+a_i)/2, (n+b_i-a_i)/2)$. Since $a_i <a_{i+1}$ and
$(a_{i+1},b_{i+1})$ is a vertex of $C_n$, it follows that
$(a_{i+1},b_{i+1})$ must actually lie inside the smaller triangle
with vertices $(a_i,b_i)$, $(a_i,n-a_i)$ and
$((n-b_i+a_i)/2, (n+b_i-a_i)/2)$.
\end{proof}

\subsection{On the difference $b_s-a_s$}
The inequalities in Proposition~\ref{ver-cond} may seem obvious,
but they play a key role in our algorithms to compute the vertices
of $C_n$. The vertex $(a_s,b_s)$ has an important
property. Let $M(n)$ denote the quantity
\begin{equation}
\label{eq:Max Diff}
                      M(n)= \max \left\{  \left\vert a-b\right\vert :1 \leq
a,b\leq
n-1 \textrm{ and } ab\equiv 1  \pmod n\right\}.
\end{equation}
An immediate consequence of Proposition~\ref{ver-cond} is that
$$b_{s}- a_{s}=M(n).$$

The quantity $M(n)$ has been studied
in~\cite{FKSY,Kh,KhShp}.
It is shown in~\cite{KhShp} that
\begin{equation}
\label{eq:MK-IS bound}
n-M(n) \ll n^{3/4 +o(1)}.
\end{equation}

On the other hand, by~\cite[Theorem 3.1]{FKSY}, for almost all $n$ $$
n-M(n) \gg n^{1/2}\( \log n \)^{\delta/2}\left ( \log\log n\right
)^{3/4}f(n), $$ where $$\delta = 1-\frac{1+\log\log 2}{\log 2}
=0.086071\ldots , $$ and $f(x)$ is any positive function tending
monotonically to zero as
$x \rightarrow \infty$. We recall that it has been proposed
in~\cite[Conjecture~4.1]{FKSY} that the above bound is quite tight:

\begin{conjecture}
\label{M-conj}
For almost all $n$
$$
n-M(n) \ll
n^{1/2}\( \log n \)^{\delta/2}\left ( \log\log n\right )^{3/4}g(n),
$$
where $g(x)$ is any function tending monotonically to $\infty$ as $x
\rightarrow \infty$.
\end{conjecture}

In support of Conjecture~\ref{M-conj} we make the following observation.
For
a fixed $\varepsilon >0$ define the set $$
\cN(\varepsilon) = \{ n \in \N \ :  \ \exists \
d|(n-1) \textrm{ such that }
n^{1/2-\varepsilon} \le d \le n^{1/2} \}.
$$
By \cite[Theorem~22]{HT} $\cN(\varepsilon)$ has positive
asymptotic density. Since
$$d\left( n - \frac{n-1}{d} \right) \equiv 1 \pmod n,$$
we see that
$$n - M(n) \le n-\left( n - \frac{n-1}{d}-d \right)=\frac{n-1}{d}+d
\ll n^{1/2 + \varepsilon},$$
for every $n$ with this property.
This immediately implies that for any $\varepsilon>0$
$$
n-M(n) \le n^{1/2+\varepsilon}
$$
for a set of $n$ of positive density, which is a weaker form
of what is assumed in Conjecture~\ref{M-conj}. In~\cite{FKSY}, one
can also find more developed heuristic arguments supporting
Conjecture~\ref{M-conj}.

We make one other remark about the vertex $(a_s,b_s)$.
Following~\cite{Ten1}, we introduce the quantities
$$\rho_1(m)=\max_{d|m,\,d\le \sqrt{m}}d \qquad \text {and}
\qquad \rho_2(m)=\min_{d|m,\,d \ge \sqrt{m}}d.$$
We note that
$$a_s=\rho_1(kn-1)
\qquad \text {and}
\qquad (n-b_s)=\rho_2(kn-1), $$
where $k$ is the integer such that $a_s(n-b_s)=kn-1$.

\subsection{Heuristic}

Our heuristic attempt to approximate $v(n)$  makes use of a
probabilistic model.
Specifically, to view the points of $n^{-1}G_n$ as being randomly
distributed in the unit square (which is supported by theoretic
results
of~\cite{CZ,GrShZa,VZ,Zha,Zhe}) and then appeal to a result of
R\'{e}nyi and Sulanke~\cite[Satz~1]{ReSu1}.
Let $\cR$ be a convex polygon in the plane with $r$ vertices and let
$P_i$, $i=1,\ldots ,n$,  be $n$ points chosen  at random in $\cR$
with uniform distribution. Let $X_n$ be the number of sides of the
convex closure of the  points $P_i$,  and let $E(X_n)$ be the
expectation of $X_n$. Then
\begin{equation} \label{eq:Expected Number}
                      E(X_n) = \frac{2}{3}r(\log n + \gamma) + c_\cR + o(1),
\end{equation}
where $\gamma= 0.577215\ldots$ is the Euler constant,
and $c_\cR$ depends on $\cR$ and is maximal when $\cR$ is a
regular $r$-gon or is affinely equivalent to a regular $r$-gon.
In particular, for the
unit square $\cR = [0,1]^2$ we have
$$
c_\cR=  - \frac{8}{3}  \log 2.
$$
More precise results are given by Buchta and Reitzner~\cite{BuRe}, but
they do not affect our arguments.

Using~\eqref{eq:Expected Number} with $r = 4$, it
is plausible to conjecture that for most $n$
\begin{equation} \label{eq:Heuristic Estimate}
v(n)   \approx  h(n),
\end{equation}
where
$$
h(n) =  \frac{8}{3}  (\log \varphi(n)  + \gamma  - \log 2 ). $$
A portion of our work has been to generate numerical data to test this
conjecture.

\section{Bounds on $v(n)$}

\subsection{Lower Bounds}

Here we give a lower bound on $v(n)$ in terms of the number of
divisors function $\tau(n)$. We begin by establishing some notation
and making a couple of pertinent observations.

For a fixed $n$, let us consider the curves $\alpha_j(n)$ and
$\beta_j(n)$ defined by \begin{eqnarray*}
\alpha_j(n): & \quad & x(n-y)=jn-1, \ 1\leq x \leq y \leq n-1,\\
\beta_j(n): & \quad &y(n-x)=jn-1, \ 1\leq y \leq x \leq n-1.
\end{eqnarray*}
A key observation used repeatedly is that for each point of $G_n$
there is a $j$ in the range $1,\ldots,\rf{n/4}$
                  such that the point lies on the curve $\alpha_j(n)$ or
$\beta_j(n)$.
We denote the region bounded by the curves $\alpha_1(n)$ and
$\beta_1(n)$ by $\cR_n$. The next figure is an illustrative example.
We note that the outermost curves are $\alpha_1(41), \beta_1(41)$.

\centerline{ \includegraphics*[scale=.45]{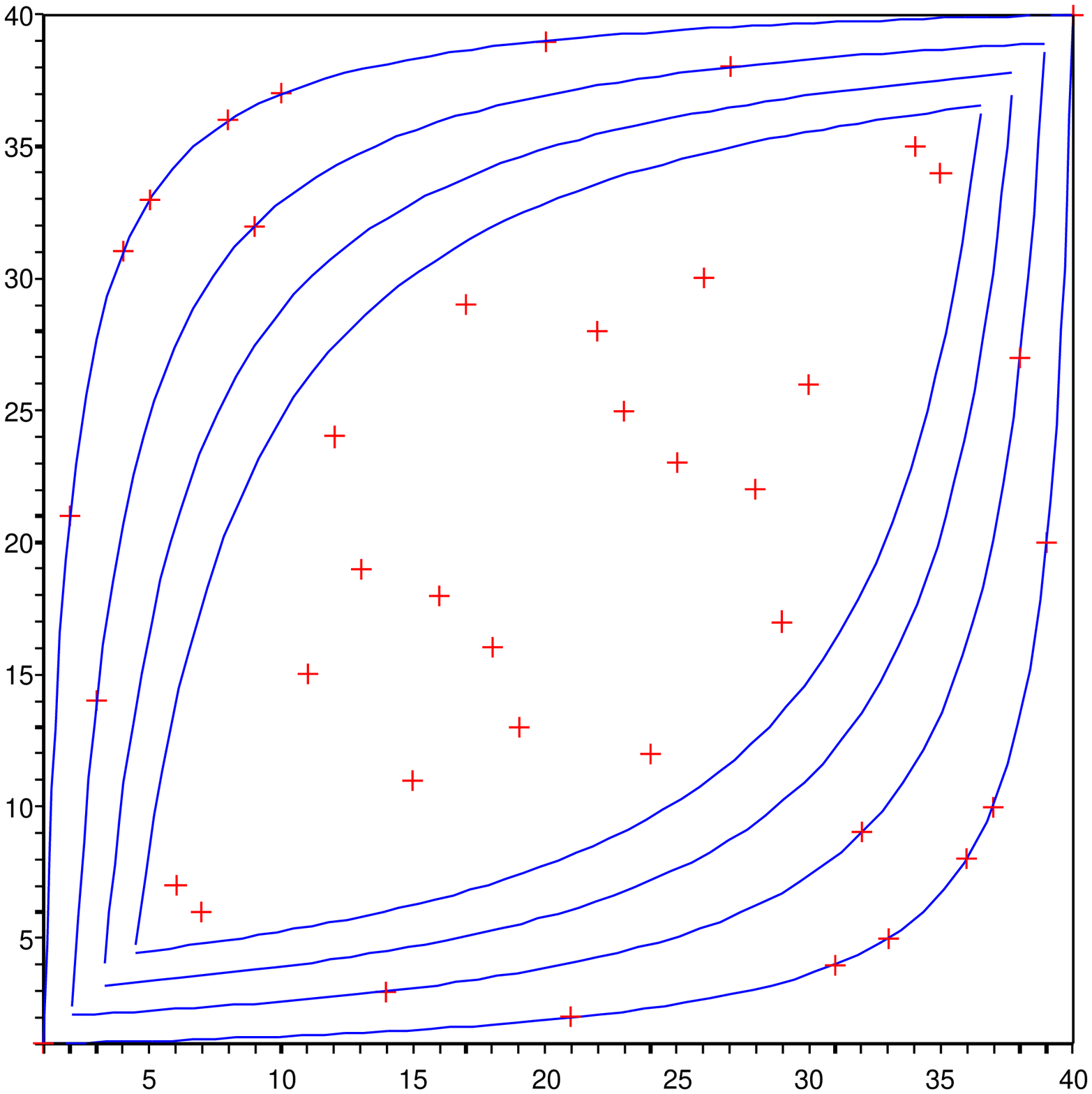} }
\centerline{{\bf Fig. 2.} The graph $G_{41}$ and the curves
$\alpha_j(41),\beta_j(41)$, $j=1,2,3,4$}
\medskip

For an integer $s\ge 1$ we denote
$$
T(s) =  \max_{i=1,\ldots,\tau(s)-1} \frac{d_{i+1}}{d_i}
$$
where $1= d_1<\ldots<d_{\tau(s)}=s$ are the positive divisors of $s$.

Clearly,
\begin{equation}
\label{eq:T and P}
T(s) \le  P(s),
\end{equation}
where $P(s)$ denotes the largest prime divisor of $s$.

Let $D_n $ be the convex closure of the points
$(d_i,n-(n-1)/d_i),(n-(n-1)/d_i,d_i)$,
for $i=1,\ldots, \tau(n-1)$. Clearly, we have the inclusions
$D_n \subseteq C_n \subseteq \cR_n$. We remark that if $n-1$ is prime,
the set $D_n$ is simply the line
segment connecting the points $(1,1)$ and $(n-1,n-1)$.

The purpose of
our next proposition is to give a criterion to determine which of the
$\alpha_j(n)$, $2 \le j \le \lceil n/4 \rceil$, lie strictly in the
interior of $D_n$, and hence strictly in the interior of $C_n$.
We denote by $\Gamma_n$ the set of boundary points $ (x,y)$ of $D_n $
such that $y \ge x$,
that is, $\Gamma_n = \{(x,y)\ : \ (x,y) \in \partial D_n,y\ge x\}$.

\begin{proposition}\label{intersec-bd}
Let $1= d_1<\ldots<d_{\tau(n-1)}=n-1$ be the positive divisors of $n-1$.
Then, for any integer $m \ge 2$,
$$\Gamma_n \cap \alpha_m(n) = \emptyset \quad \Leftrightarrow \quad
\frac{d_{i+1}}{d_i} + \frac{d_i}{d_{i+1}} < 4m-2 +\frac{4(m-1)}{n-1},
\ i=1,\ldots,\tau(n-1)-1.
$$
\end{proposition}

\begin{proof}
This is a routine computation and so we only sketch an outline.
The polygonal curve $\Gamma_n$ is the union of line segments
$$ L_i:
\quad (1-t)(d_i,n-(n-1)/d_i)+t(d_{i+1},n-(n-1)/d_{i+1}),\, 0\le t \le
1,$$
with $i=1,\ldots,(\tau(n-1)-1)$. Now $L_i \cap \alpha_m(n) = \emptyset$ if
and only if the quadratic equation in $t$ $$
(d_{i+1}-d_i)\(\frac{n-1}{d_{i+1}}-\frac{n-1}{d_i}\)t^2-
(d_{i+1}-d_i)\(\frac{n-1}{d_{i+1}}-\frac{n-1}{d_i}\)t+(1-m)n=0
$$
has no real solutions.
\end{proof}

A useful consequence of Proposition~\ref{intersec-bd} is that if
\begin{equation}
\label{eq:m and T(n-1)}
m \ge \fl{ \frac{T(n-1) + 3}{4}},
\end{equation}
with $m \in \Z$ and $m \ge 2$, then
$\Gamma_n \cap \alpha_m(n) = \emptyset.$

\begin{theorem}
\label{thm:Lower Bound}
For all $n \ge 2$,
$$
v(n)\geq 2\( \tau(n-1)-1 \),
$$
and for sufficiently large $x$,
$$
\# \left\{n \le x\ : \  v(n) = 2\( \tau(n-1)-1 \)\right\} \gg \frac{x}{\log x}.
$$
\end{theorem}

\begin{proof}
Since $C_n \subseteq \cR_n$, any $(x,y) \in G_n \cap \left
(\alpha_1(n) \cup \beta_1(n) \)$
is a vertex of $C_n$, and either $x$ or $y$ is a divisor of $(n-1)$.
Therefore, $v(n)\geq 2\(\tau(n-1)-1\)$.

By~\eqref{eq:m and T(n-1)} we
have $\Gamma_n \cap \alpha_2(n) = \emptyset$ for every $n$ with $T(n-1) \le 5$.
        Consequently, for such $n$, all of the
vertices of $C_{n}$ lie on $\alpha_1(n) \cup \beta_1(n)$ and thus
$v(n) = 2\( \tau(n-1)-1 \)$.
On the other hand, by~\cite[Theorem~1]{Saias}, we know that
for any fixed $t$ and sufficiently large $x$,
$$
\# \left\{n \le x\ : \  T(n-1) \le t \right\} \asymp
\frac{x\log t}{\log x}
$$
Applying this result with $t=5$ we conclude the proof.
\end{proof}

It is easy to construct explicit examples of $n$ with
$v(n) = 2\( \tau(n-1)-1 \)$.
For instance it follows from~\eqref{eq:T and P} and~\eqref{eq:m and T(n-1)}
that this holds for  $n=2^r3^s5^t+1$,
where $r,s,t$ are non-negative integers.

Since for any $\delta>0$ we have
               $$
\limsup_{k\rightarrow\infty} \tau(k) 2^{-(1-\delta)\log k/ \log \log
k} = \infty
$$
(see~\cite[Theorem~317]{HarWr}), the same holds true for $v(n)$, and so
we can
infer that the heuristic estimate~\eqref{eq:Heuristic Estimate} is
sometimes
exponentially smaller than $v(n)$.

\begin{corollary}
\label{cor:Large v}
For any $\delta>0$
$$
\limsup_{n\rightarrow \infty}  v(n) 2^{-3/8(1-\delta)h(n)/ \log h(n)}
=  \infty.
$$
\end{corollary}

We have that $v(n) \ge 2(\tau(n-1)-1)$, and it is natural to ask when
does one have strict inequality. Our next result gives a partial
answer to this question. Specifically, we  exhibit a set of
positive density for which we have strict inequality.  Furthermore, if
we assume Conjecture~\ref{M-conj} then we have strict inequality for
almost all $n$.

\begin{theorem}
\label{thm:Strict Ineq}
The strict inequality
$$ v(n) > 2(\tau(n-1)-1)$$
holds
\begin{itemize}
\item[i.]  for a set of $n$ of positive density.
\item[ii.] for almost all $n$,
provided that for almost all $n$ we have
$n-M(n) \le n^{1/2 + o(1)}$.
\end{itemize}
\end{theorem}

\begin{proof}
{\emph {i.}} Let
$$\cE(x)= \{n\le x \ : \ v(n)=2(\tau(n-1)-1) \}, $$
and $$\cI(x)= \left \{ n \le x\ : \   a_s\(n-b_s\)= n-1 \right \}. $$
It is important to note that the values of $s,\ a_s$ and $b_s$ all
depend on $n$. We remind the reader of the following properties of the
point $(a_s,b_s)$ used in the proof below.
It is the highest vertex of $C_n$ that lies on or below the line
$x+y=n$; $M(n) = b_s-a_s$ and $a_s \le n-b_s$. Clearly, $\cE(x)
\subseteq \cI(x)$.

The set of positive density we have in mind is
$$\cA(x) = \{n \le x \ : \ \exists p\ \text{prime with}\ p|(n-1)\
\text{and}\ p
\ge x^{0.76} \}.$$
Using {\it  Mertens's formula\/},
(see~\cite[Theorem~427]{HarWr}), we get that
$$ \# \cA(x) = \sum_{x^{0.76} \le p \le x}\fl{\frac{x-1}{p}}
\sim (\log(25/19))x.$$

Since $\cE(x) \subseteq \cI(x)$, in order to prove
$$ \lim_{x \rightarrow \infty}\frac{\#(\cA(x)\cap \cE(x))}{x} =0 $$
it is enough to prove that
$$ \lim_{x \rightarrow \infty}\frac{\#(\cA(x)\cap \cI(x))}{x} =0. $$

We now write $\cI(x)$ as the disjoint union of the two sets $\cI_1(x)$,
$\cI_2(x)$,
where
\begin{eqnarray*}
\cI_1(x) &= &\left \{ n \in \cI(x) \ : \  n-b_s \le
x^{0.24} \right \}, \\
\cI_2(x) &= & \left \{ n \in \cI(x) \ : \ x^{0.24} < n-b_s < x^{0.76}
\right \}.
\end{eqnarray*}
The exponent values, 0.24 and 0.76, come from the asymptotic $n-M(n)
\le n^{3/4+o(1)}$ that we mentioned earlier. Since $\# \cI_1(x) \le
x^{0.48}$ and for $x$ large $\cA(x)\cap \cI_2(x) = \emptyset$, it
follows
that for large $x$
$$\#(\cA(x)\cap \cI(x))=\#(\cA(x)\cap \cI_1(x))+\#(\cA(x)\cap
\cI_2(x)) \le \#\cI_1(x) =o(x).$$

{\emph {ii.}} We now prove the following conditional statement. If
for almost all $n$,
$n-M(n) \le n^{1/2} g(n)$  with some function $g(n)=n^{o(1)}$, then
$\#\cI(x)=o(x).$

Without loss of generality we may assume that $g(n)$ is monotonically
increasing. This time we write $\cI(x)$ as the disjoint union of three
sets, $\cJ_1(x), \cJ_2(x)$ and $\cJ_3(x)$ where \begin{eqnarray*}
\cJ_1(x) &= &\left \{ n \in \cI(x) \ : \  n-b_s \le
\frac{\sqrt{x}}{g(x)} \right \}, \\
\cJ_2(x) &= & \left \{ n \in \cI(x) \ : \
\frac{\sqrt{x}}{g(x)} < n-b_s \le \sqrt{x}g(x) \right \},\\
                     \cJ_3(x) &= & \left \{ n \in \cI(x) \ : \  \sqrt{x}g(x) <
n-b_s < x^{0.76}
            \right \}.
\end{eqnarray*}

Now $\# \cJ_1(x) \le xg(x)^{-2} =o(x)$, and by our assumption
we also have $\# \cJ_3(x) =o(x)$. So to conclude we
need  to show that $\# \cJ_2(x)=o(x)$. This follows by the following
observation. Let
$$H(x,y,z)= \{ n \le x\ : \ \exists d|n \textrm{ with }y<d \le z \}.$$
Then $$\# \cJ_2(x) \le H\(x,\sqrt{x}/g(x), \sqrt{x}g(x) \),$$
and by~\cite[Theorem~1]{F},
$$
H\(x,\sqrt{x}/g(x), \sqrt{x}g(x) \) =o (x)
$$
which concludes the proof.
\end{proof}

We remark that the assumption of Theorem~\ref{thm:Strict Ineq}~(ii)
is weaker than Conjecture~\ref{M-conj}. The bound of 
Conjecture~\ref{M-conj} probably holds for almost all
primes. This would then imply that $$ v(p)> 2\( \tau(p-1)-1 \) $$ for
almost all primes $p$.
On the other hand, it is reasonable to expect that there are
infinitely many primes of the form $n=2^r3^s5^t+1$ (in fact even of
the form $p = 3\cdot2^r +1$), and therefore equality would occur
infinitely often, as well.
We conclude this section by proving that
$v(n)$ can be substantially larger than $\tau(n-1)$.

\begin{theorem}
\label{thm:v vs tau}
There is an infinite sequence of integers $n_j$
with
$$
v(n_j) \ge \exp\(\(\frac{2\log 2 }{11} + o(1)\) \frac{\log n_j}{\log 
\log n_j}\)
\qquad \text{and}\qquad
\tau(n_j-1) =2.
$$
\end{theorem}

\begin{proof}
Let $n$ be a shifted prime,
that is, $n=p+1$, where $p$ is prime.
We first show that for such
integers,
$$
v(n) = v(p+1)\ge 2(\tau(2p+1)-3).
$$
Let $\ell$ be the line through $(1,1)$ which is tangent to
$\alpha_2(n)$. Since $(1,1)$ and $(p,p)$ are the only points
of $G_n$ on $\alpha_1(n)$, all of
the points of $G_n$ lie on or below $\ell$.
A straightforward calculation shows that $\ell$  meets
$\alpha_2(n)$ at the point $(x,y)$  where the $x$-coordinate is
$$
x = \frac{1}{1-((p+1)/(2p+1))^{1/2}} < 4.
$$
Hence every divisor $d$ of $2p+1$, with
$3 < d < (2p+1)/3$, gives rise to a vertex on $\alpha_2(n)$.
Consequently the number of vertices on $\alpha_2(n)$ is at least 
$\tau(2p+1)-4$.
By symmetry there are
an equal number of vertices on $\beta_2(n)$, and since
$(1,1)$ and $(p,p)$ are also vertices of
$C_n$, we obtain the desired inequality.

We now let $Q_j$ denote the product of first $j$ odd primes and
set $p_j$ to be the smallest prime satisfying the congruence
$2p_j \equiv -1 \pmod {Q_j}$. By the Prime Number Theorem $\log Q_j 
\sim j \log j$, and by Heath-Brown's~\cite{He-Br} version of Linnik's 
theorem we have
$p_j < cQ_j^{11/2}$, for an absolute constant $c \ge 1$.
On combining $p_j < cQ_j^{11/2}$ with the asymptotic $\log Q_j \sim j \log j$
we obtain
$$
\tau(2 p_j +1) \ge \tau(Q_j) = 2^j \ge
\exp\(\(\frac{2\log 2 }{11} + o(1)\) \frac{\log p_j}{\log \log p_j}\).
$$
Setting $n_j = p_j +1$ we conclude the proof.
\end{proof}

In particular, we see from Theorem~\ref{thm:v vs tau} that
$$
\limsup_{n \rightarrow \infty}\frac{\log v(n)}{\log \tau(n-1)} = \infty.
$$
Furthermore we can replace the terms  $\log v(n)$ and  $\log 
\tau(n-1)$ by the $k$-fold iteration of the logarithm for any $k \in 
\N$. Unfortunately, we do not see any approaches to the following.

\begin{conjecture} We have
$$
\liminf_{n \rightarrow \infty}v(n) = \infty.
$$
\end{conjecture}

\subsection{Upper Bounds}

\begin{theorem}
\label{thm:Upper Bound-1}
For $n \to \infty$,
$$
v(n)\le n^{3/4 + o(1)}.
$$
\end{theorem}

\begin{proof}  In Section~\ref{sec:Triangle Points},
we labelled the highest vertex of $C_n$ in the triangle $\cT_n$ by
$(a_s,b_s)$.
Trivially, $s \le a_s$ and $a_s \le n-b_s$. Hence
$$v(n) \le 4s+2 \le 4a_s+2 \le 2(n-b_s+a_s+1) = 2(n-M(n)+1),$$
and the bound~\eqref{eq:MK-IS bound} concludes the proof.
\end{proof}

Most certainly the bound of Theorem~\ref{thm:Upper Bound-1} is not
tight. If we assume Conjecture~\ref{M-conj}, then $$ v(n)\le n^{1/2 +
o(1)} $$ for almost all $n$. This still seems too high and the actual
order of $v(n)$ is almost certainly much smaller.
A different upper bound for $v(n)$ can be derived from~\eqref{eq:m and T(n-1)}.
For integers $n$ where $n-1$
has only small prime factors, this upper bound is significantly
better than Theorem~\ref{thm:Upper Bound-1}.

\begin{theorem}
\label{thm:Upper Bound-2}
For $n \to \infty$,
$$
v(n)\le T(n-1) n^{o(1)}.
$$
\end{theorem}

\begin{proof} From~\eqref{eq:m and T(n-1)} we see that only points
from the curves $\alpha_j(n)$ and $\beta_j(n)$ where,
                 $$
j \le m_n =\fl{ \frac{T(n-1) + 3}{4}},
$$
contribute to $v(n)$.
Since every curve $\alpha_j(n)$, $\beta_j(n)$ contains at most
$\tau(jn-1)$ points of $G_n$ we derive $$
v(n) \le \sum_{j =1}^{m_n} 2 \tau(jn-1).
$$
               We conclude by invoking the asymptotic
inequality $\tau(r) \ll r^{o(1)}$, see~\cite[Theorem~315]{HarWr}.
\end{proof}

\section{Computing $C_n$}

\subsection{Systematic search algorithm}

We now describe a deterministic algorithm to construct the vertices
of $C_n$ that lie in the triangle $\cT_{n}$.  It is a variant  of the
famous algorithm of Graham~\cite{Grah} known as \textsc{Graham Scan}.
The main virtue of our algorithm, as  opposed to using some other
convex closure algorithms,  is that we do not need to generate and
store all of the points of $G_n$ before determining the convex
closure.  Instead, we generate the points one by one, discard most of
them along the way, and halt in a reasonable amount of time.

\begin{algorithm}  \label{alg:Search Alg}
\begin{enumerate}
\item Set $a_{0}:=1;b_{0}:=1$.
\item For $i =0,1, \ldots $:
\begin{enumerate}
\item Set $a_{i+1}:=$ to be the smallest integer $a \in \Z_{n}^{\ast}$
satisfying the inequalities $$ a_{i}<a \le \frac{n+a_{i}-b_{i}}{2}
\quad \text{and} \quad b_{i}-a_{i} < a^{-1}-a.
$$
If either of the above conditions cannot be met the algorithm
terminates.
\item Set $b_{i+1}:=a^{-1}$.
\item Convexity check: %Obtain the convex closure:
\begin{enumerate}
\item If $i =1$ goto Step~2(a).
\item If $i \geq 2$ and the angle between the points $\(
a_{i-1},b_{i-1}\),\(  a_{i},b_{i}\)$ and $\(a_{i+1},b_{i+1}\)$  is
reflex then return to Step~2(a), otherwise discard the point $\(
a_{i},b_{i}\)$ and set $$a_{i}:=a_{i+1}, \quad b_{i}:=b_{i+1}, \quad
i:
= i-1$$ and return to Step~2(c).
\end{enumerate}
\end{enumerate}
\end{enumerate}
\end{algorithm}

              We note that the inequalities in Step~2a are
motivated by Proposition~\ref{ver-cond}. Clearly,
Algorithm~\ref{alg:Search Alg} is deterministic and it
immediately follows from~\eqref{eq:MK-IS bound} that its complexity is
$O(n^{3/4 + o(1)})$.

\subsection{Factorisation based algorithm}

The observation that the points in $G_n\cap \alpha_1(n)$ are vertices
of $C_n$ combined with~\eqref{eq:m and T(n-1)} allows us to devise a
variation on
Algorithm~\ref{alg:Search Alg}.
The idea is to first use factorisation to create a smaller input set
and then run the algorithm.

Let $\cP_n$ be the polygonal region with vertices
\begin{align*} &(1,n-1),(1,1),
\(d_1,n-(n-1)/d_1\),\ldots,\(d_k,n-(n-1)/d_k\),\\
&\(\((n-1)/d_k+d_k\)/2,n-\((n-1)/d_k+d_k\)/2\right
), (\sqrt{n-1},n-\sqrt{n-1}),
\end{align*}
where  $1=d_0<d_1<\ldots <d_k $ are the factors of $n-1$ which are
less than or equal to
$\sqrt{n-1}$.
Since the vertices of $C_n$ can only lie on the
curves $\alpha_j(n)$, $\beta_j(n)$  where
$$
j \le m_n = \fl{\frac{T(n -1)+3}{4}},
$$
we need only determine which of the points
of the union $$ U_n = \bigcup_{j=1}^{m_n} S_{j,n}, $$ are vertices of
$C_n$, where $S_{j,n}=
\alpha_j(n)\cap G_n \cap \cP_n$.
It is useful to keep in mind that
$$
\# U_n \le \sum_{j =1}^{m_n} \# S_{j,n} \le \sum_{j =1}^{m_n}
\tau(jn-1) = m_n n^{o(1)},
$$
see~\cite[Theorem~315]{HarWr}. We now apply the following algorithm.

\begin{algorithm}  \label{alg:Fact Alg} \ \begin{enumerate}
\item %Step 1:
Factorization: \label{stp:fact}
\begin{enumerate}
\item Find
all of the factors $1=d_0<d_1<\ldots <d_k\le  \sqrt{n-1}$ of $n-1$.
\item Set $S_1 := \{(1,1),
\(d_1,n-(n-1)/d_1\),\ldots,\(d_k,n-(n-1)/d_k\)\}$.
\item  Compute $t:=  T(n -1)$.
\item Set  $m_n:=  \fl{(t+3)/4}$. \label{stp:bound-on-levels}
\item For $j=2, \ldots, m_n$, factor $jn -1$ and construct the set
$S_{j,n}$.
\item Set $U_n := \cup_{j=1}^{m_n} S_{j,n}$.
\end{enumerate}
\item %Step 2:
Determining the vertices: \label{stp:ver-det}
\begin{enumerate}
\item Order the points of $U_n$ by increasing first co-ordinate.
\item Apply the appropriate versions of Steps 2a and 2c of
Algorithm~\ref{alg:Search Alg} to the elements of $U_n$.
\end{enumerate}
\end{enumerate}
\end{algorithm}

The
complexity of Algorithm~\ref{alg:Fact Alg}
depends on the type of algorithm we use for the factorisation step. 
If we use any
subexponential probabilistic factorisation algorithm which runs in
time $n^{o(1)}$, (see~\cite[Chapter 6]{CrPom}), then  the complexity
of Step~\ref{stp:fact} of Algorithm~\ref{alg:Fact Alg} is at most
$$
              \# U_n  n^{o(1)}   =   m_n n^{o(1)} .
$$
Furthermore, the complexity of Step~\ref{stp:ver-det} of
Algorithm~\ref{alg:Fact Alg}
is of the same form as well.
So the
overall complexity of  Algorithm~\ref{alg:Fact Alg} is at most
$$
              m_nn^{o(1)}  =   T(n-1) n^{o(1)} .
$$

This is lower than that of Algorithm~\ref{alg:Search Alg} if $T(n-1)
\le n^{3/4}$.
For any fixed $\lambda \ge 0$ the proportion of the positive
integers
$k$ with $T(k) \le k^{\lambda}$ is given by a certain continuous function
$\psi(\lambda)>0$, see~\cite{Ten2}. Using~\cite[Corollary~A]{Saias} we conclude
that
$$
\psi(3/4) = \int_0^{7/8} \rho\(\frac{1}{x} - 1\) \frac{d\,x}{x}
= \int_{1/7}^\infty \rho\(y\) \frac{d\,y}{1+ y} =0.866468\ldots
$$
where $\rho(u)$ is the
{\it Dickman  function\/}, see~\cite{Dic}
or~\cite[Section~III.5.4]{Ten3}. Thus the proportion of the positive
integers $n$ with $T(n-1) \le n^{3/4}$ is $\psi(3/4) = 0.866468\ldots$.
(The bound in Step~\ref{stp:bound-on-levels}
of Algorithm~\ref{alg:Fact Alg} is certainly not
tight. It can probably be replaced by a bound of order $n^{o(1)}$ or
even possibly a power of $\log n$, but unfortunately
we have not been able to prove such a result.)

On the other hand, if we use a deterministic factoring algorithm
in Step~1,
then Algorithm~\ref{alg:Fact Alg} is of complexity at most
$$ m_n (m_nn)^{1/4 + o(1)} =
T(n-1)^{5/4}n^{1/4 + o(1)} $$
unconditionally, and of complexity  at most
$$
m_n (m_nn)^{1/5 + o(1)} = T(n-1)^{6/5}n^{1/5 + o(1)}
$$
under the {\it Extended Riemann Hypothesis\/}, see~\cite[Section~6.3]{CrPom}.
Accordingly,  this is better than Algorithm~\ref{alg:Search Alg} for
$T(n-1) < n^{2/5}$ and $T(n-1) < n^{11/24}$ respectively.
The corresponding proportions of the positive integers, $n$, satisfying
these inequalities are $\psi(2/5)$ and $\psi(11/24)$.
Since~\cite[Corollary~A]{Saias} expresses both $\psi(2/5)$ and
$\psi(11/24)$ as double integrals, it is easier to compute
$\psi(3/4)$ than either of these two values.

\section{Computational Results}

\subsection{Expected value of $V(N)$}

Let
$$
\eta =  \sum_{p}  \frac{\log(1-1/p)}{p} = -0.580058\ldots\,,
$$
where  the sum runs over all prime numbers $p$. Surprisingly
enough, this quantity has already appeared in various, seemingly
unrelated number theoretic questions, see~\cite[page~122]{Fin}.

\begin{proposition}
\label{eq:Sum Euler} We have,
$$
\frac{1}{N}  \sum_{n=1}^N \log \varphi(n) = \log N  + \eta -
1+ O\left(\frac{\log \log N}{N}\right).
$$
\end{proposition}

\begin{proof}
Obviously,
$$
\frac{1}{N}  \sum_{n=1}^N \log \varphi(n) = \frac{1}{N}  \sum_{n=1}^N
\log n  + \frac{1}{N}  \sum_{n=1}^N \sum_{p|n} \log(1-1/p), $$ where
the last sum is taken over prime divisors $p |n$.
The first sum on the right-hand side is $\log N - 1 + o(1)$ by {\it
Stirling's formula\/}.
By changing the order of summation in the second sum, we derive
\begin{eqnarray*}
                     \frac{1}{N}  \sum_{n=1}^N  \sum_{p|n} \log(1-1/p) & = &
\frac{1}{N} \sum_{p\le N}  \log(1-1/p) \sum_{\substack{n \le N\\p|n}}
1\\ & = & \frac{1}{N}  \sum_{p\le N}  \log(1-1/p) \(\frac{N}{p} +
O(1)\) \\
& = &     \sum_{p\le N}  \frac{\log(1-1/p)}{p} + O\(\frac{1}{N}
\sum_{p\le N} \frac{1}{p}\)\\
& = &   \sum_{p\le N}  \frac{\log(1-1/p)}{p} +  O\left(\frac{\log
\log N}{N}\right),
\end{eqnarray*}
where the last step follows by {\it Mertens's formula\/},
see~\cite[Theorem~427]{HarWr}.
                    Observing that
$$
                          \sum_{p\le N} \frac{\log(1-1/p)}{p} =  \eta 
- \sum_{p >
N} \frac{\log(1-1/p)}{p}
                          = \eta + O\left(\frac{1}{N}\right),
$$
we conclude our proof.
\end{proof}

Combining heuristic~\eqref{eq:Heuristic Estimate} with
Proposition~\ref{eq:Sum Euler} for the average $V(N)$,
we get the heuristic $V(N)\sim H(N),$ where
$$
H(N) =  \frac{8}{3}  (\log N + \gamma +
\eta  -1 - \log 2) \approx 2.66666 \cdot\log N - 4.52264.
$$

In Figure~3 we
compare the graph of $V(N)$, $H(N)$
and the least squares approximation
\begin{equation}
\label{eq:LS approx}
L(N)=3.551166\cdot \log N - 9.610899
\end{equation}
to $V(N)$, where $N$ ranges
over the interval $[2,5770001]$. The values of $V(N)$ are represented
by diamonds along the graph of $L(N)$,
while $H(N)$ is the lower curve.

\centerline{ \includegraphics*[scale=.5]{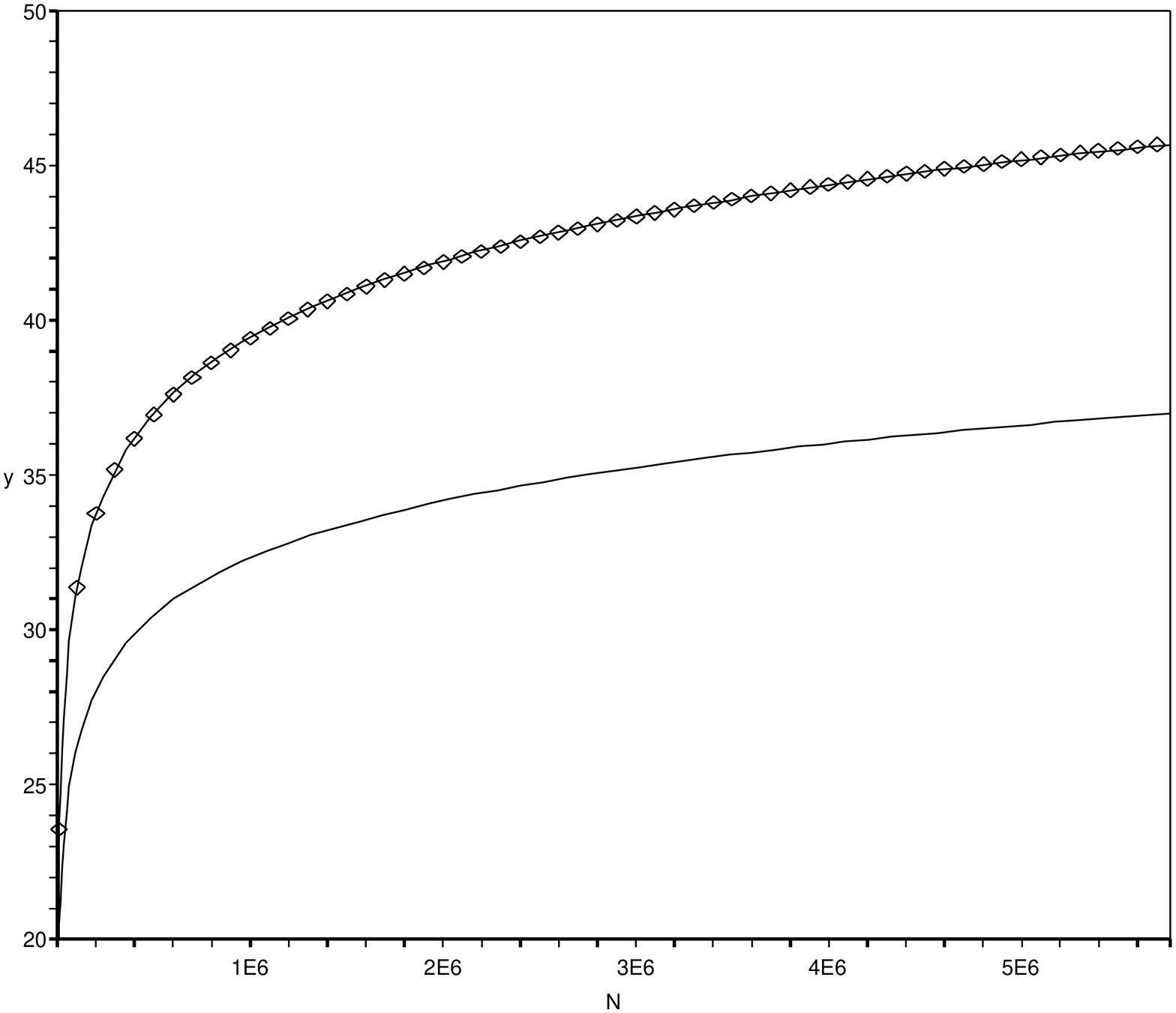} }

\centerline{{\bf Fig. 3.} $ V(N)$, $H(N)$, and $L(N)$ for
$2 \le N \le 5770001$}
\medskip

We see that although $V(N)$ behaves  like a logarithmic function
and thus resembles $H(N)$, they clearly deviate.
This deviation seems to be of regular nature and suggests
that there should be a natural explanation for this behaviour
of $V(N)$. In an attempt to understand this we computed $v(n),h(n)$
and $\tau(n-1)$ for 50000 random integers in the interval $[10^6,
10^8]$, and did some comparisons. We present the individual data in
the histograms in
Figures~4 and~5, and the comparisons in Figures~6, 7, 8, 9 and~11.
In several histograms the extreme values on the right are not
visible. Hence, for visual clarity we have truncated them on the
right. Under each histogram we state in the caption the minimum
value, the maximum value and the number of values that are not shown.

\bigskip
\centerline{ \includegraphics*[scale=.5]{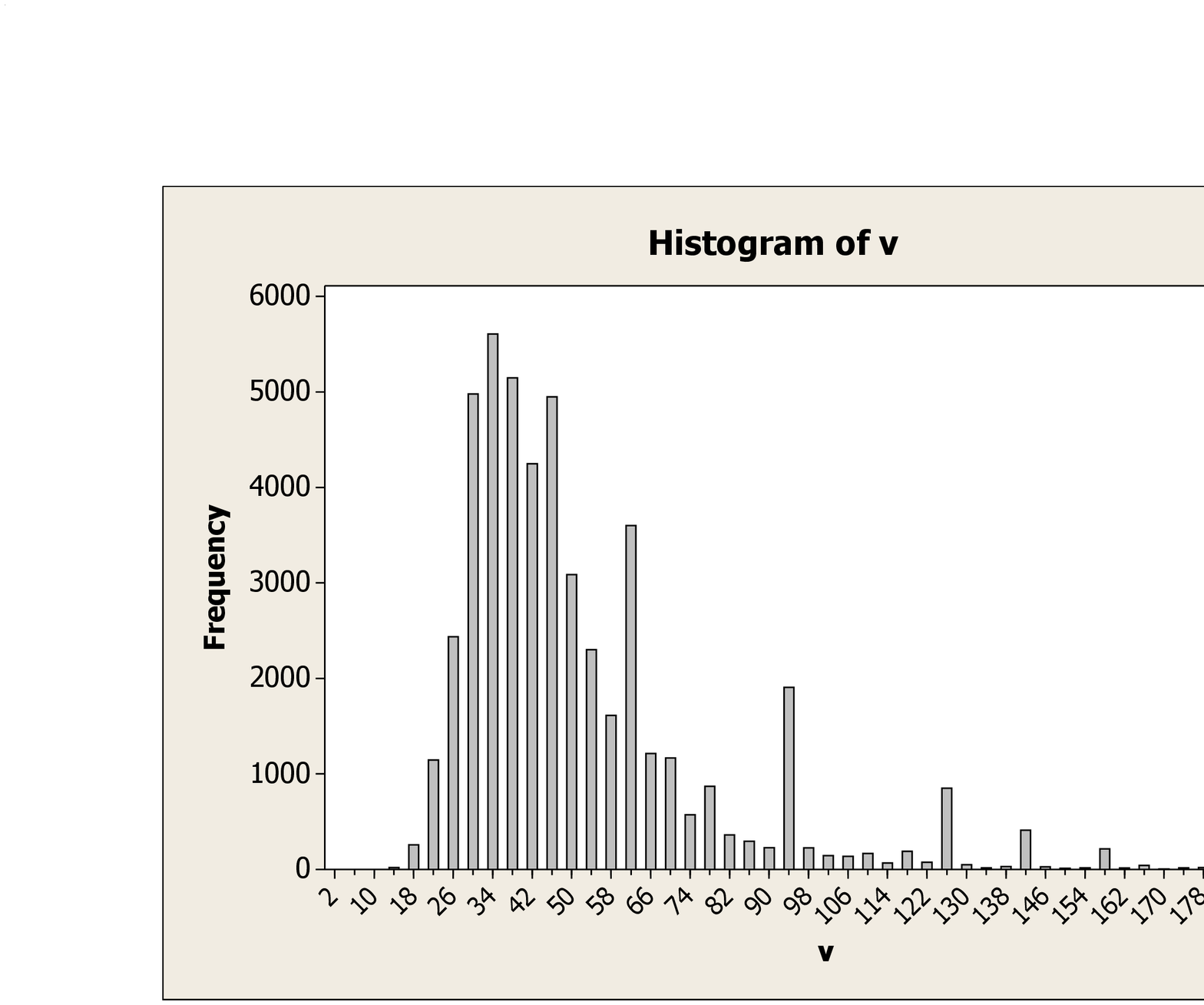} }

\centerline{{\bf Fig. 4.} Frequency histogram of $v(n)$}
\centerline{$\min=14,\ \max=766$ \quad (645 values omitted)}
\medskip

\centerline{ \includegraphics*[scale=.5]{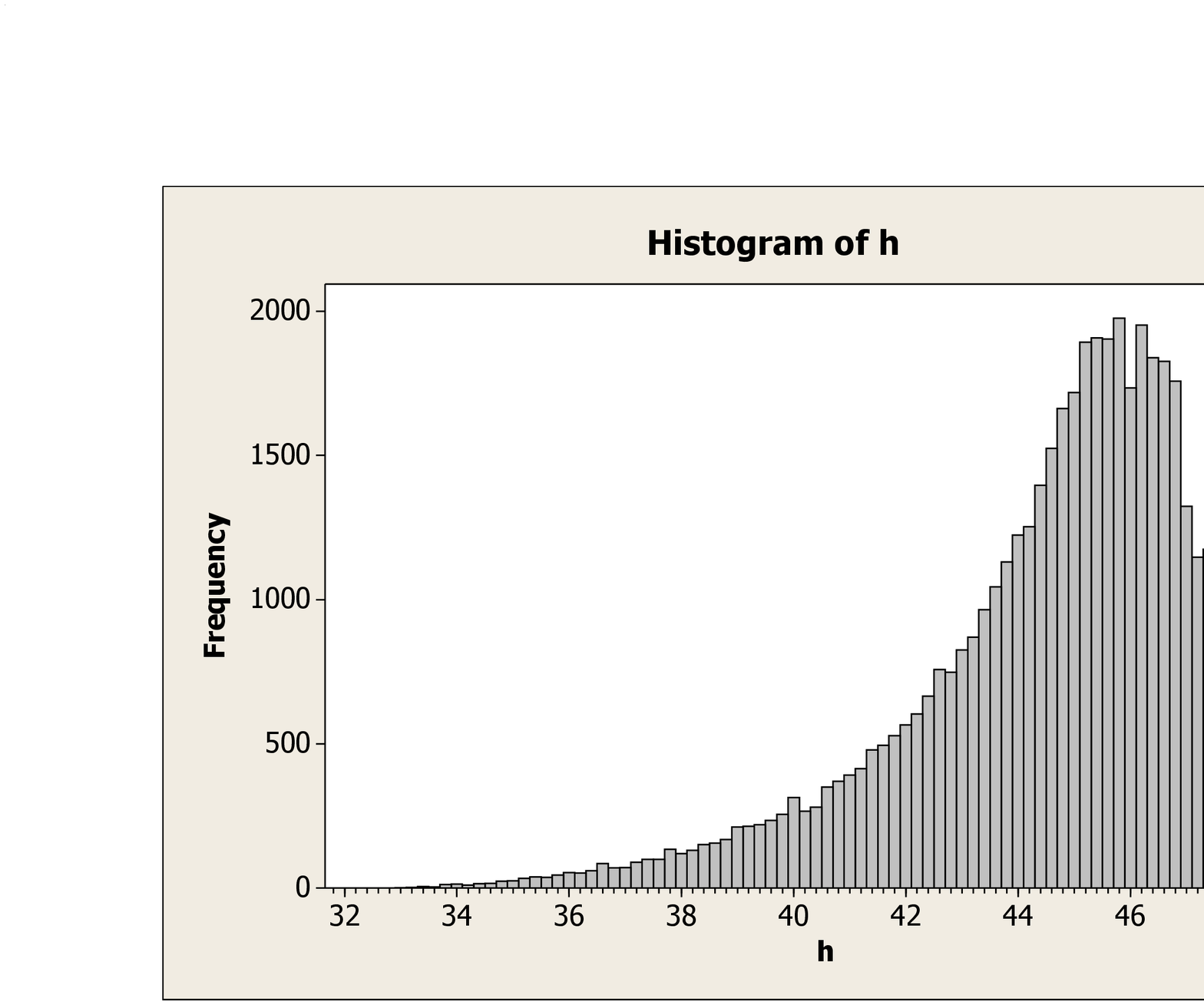} }

\centerline{{\bf Fig. 5.} Frequency histogram of $h(n)$}
\centerline{$\min=33.01,\ \max=48.81$}
\medskip

\centerline{ \includegraphics*[scale=.5]{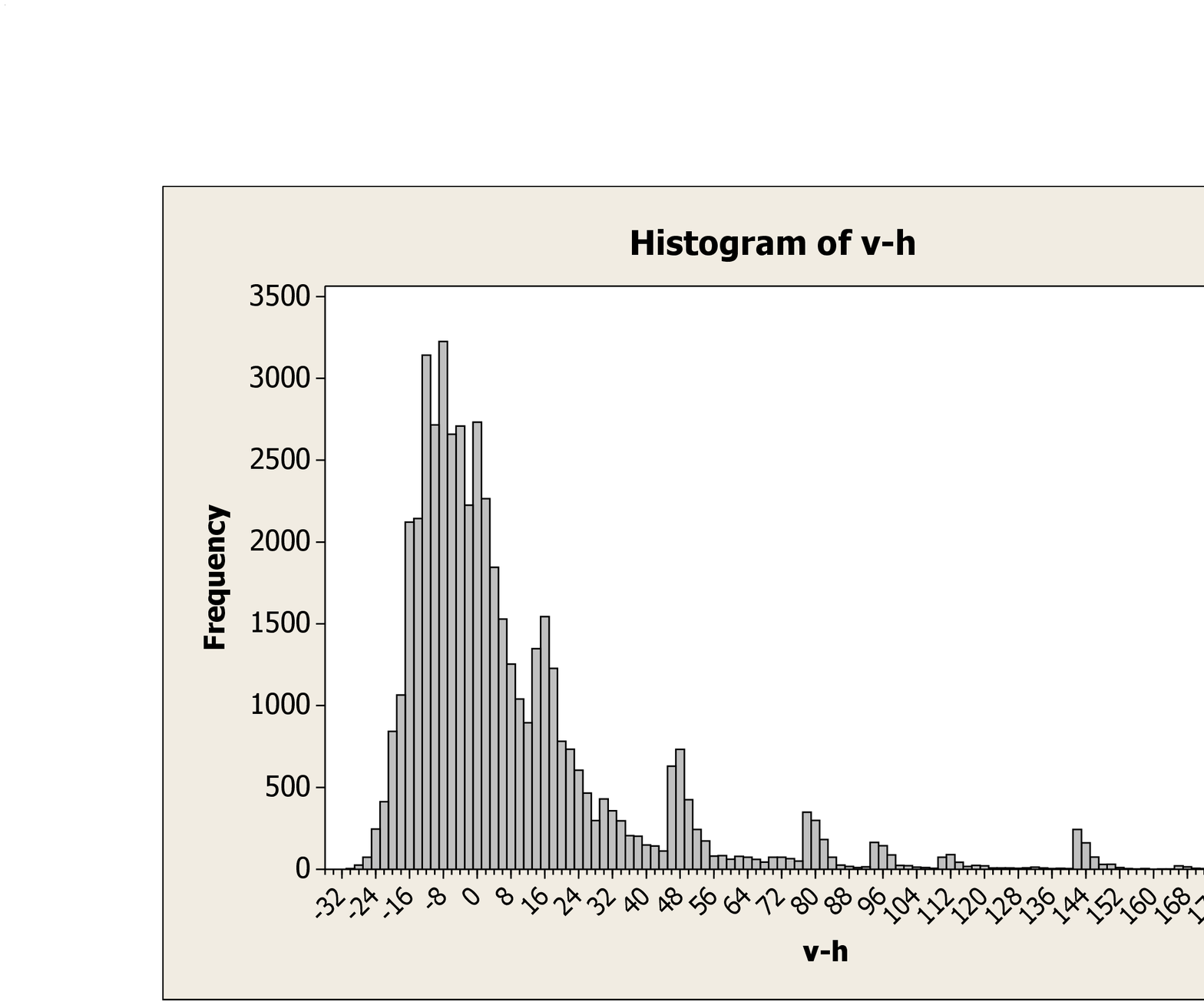} }

\centerline{{\bf Fig. 6.} Frequency histogram of $(v-h)$}
\centerline{$\min=-29.93,\ \max=714.41$ \quad (458 values omitted)}
\medskip

\centerline{ \includegraphics*[scale=.5]{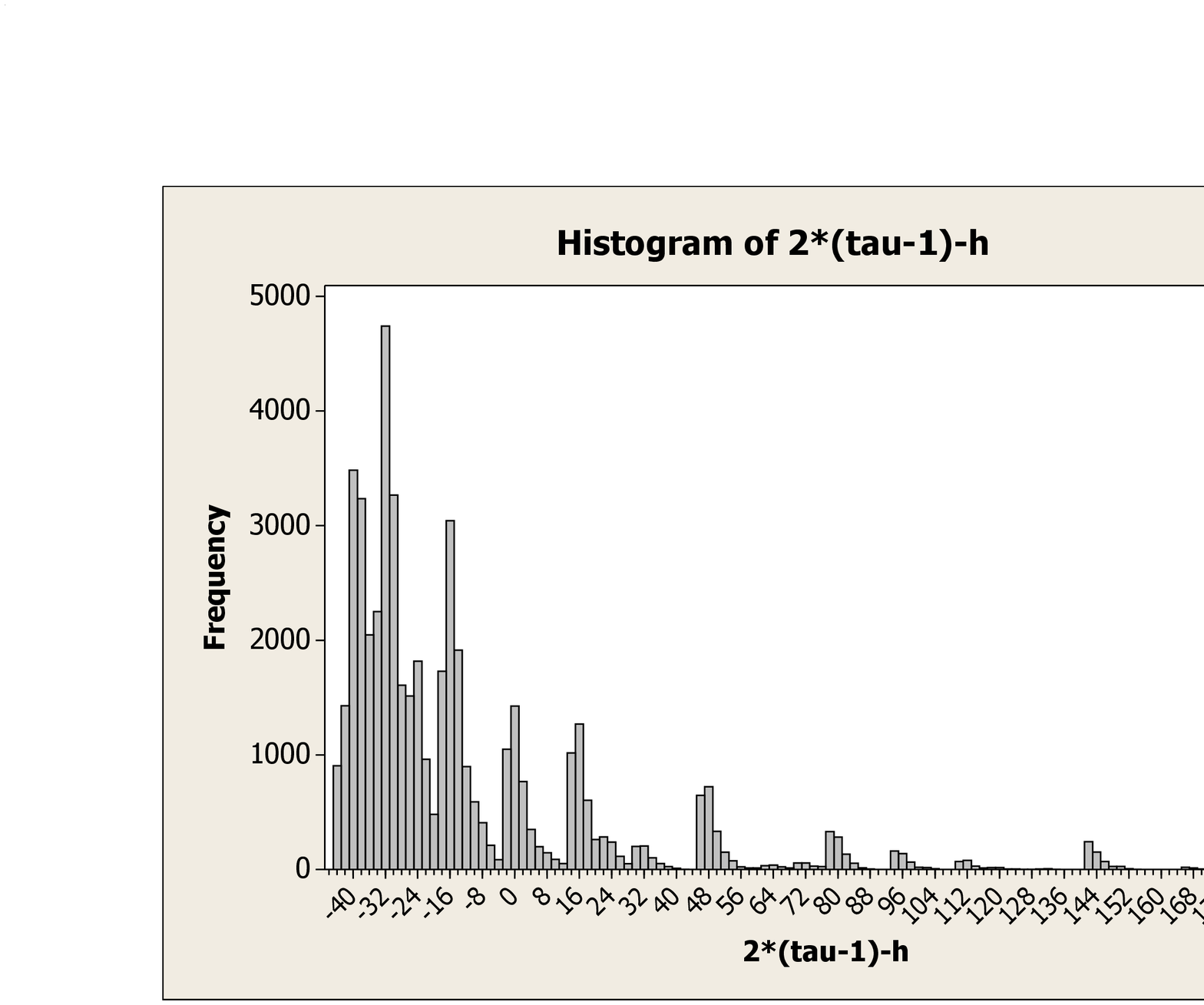} }

\centerline{{\bf Fig. 7.} Frequency histogram of $2(\tau(n-1)-1)-h(n)$}
\centerline{$\min=-44.96,\ \max=714.41$ \quad (443 values omitted)}
\medskip

\centerline{ \includegraphics*[scale=.5]{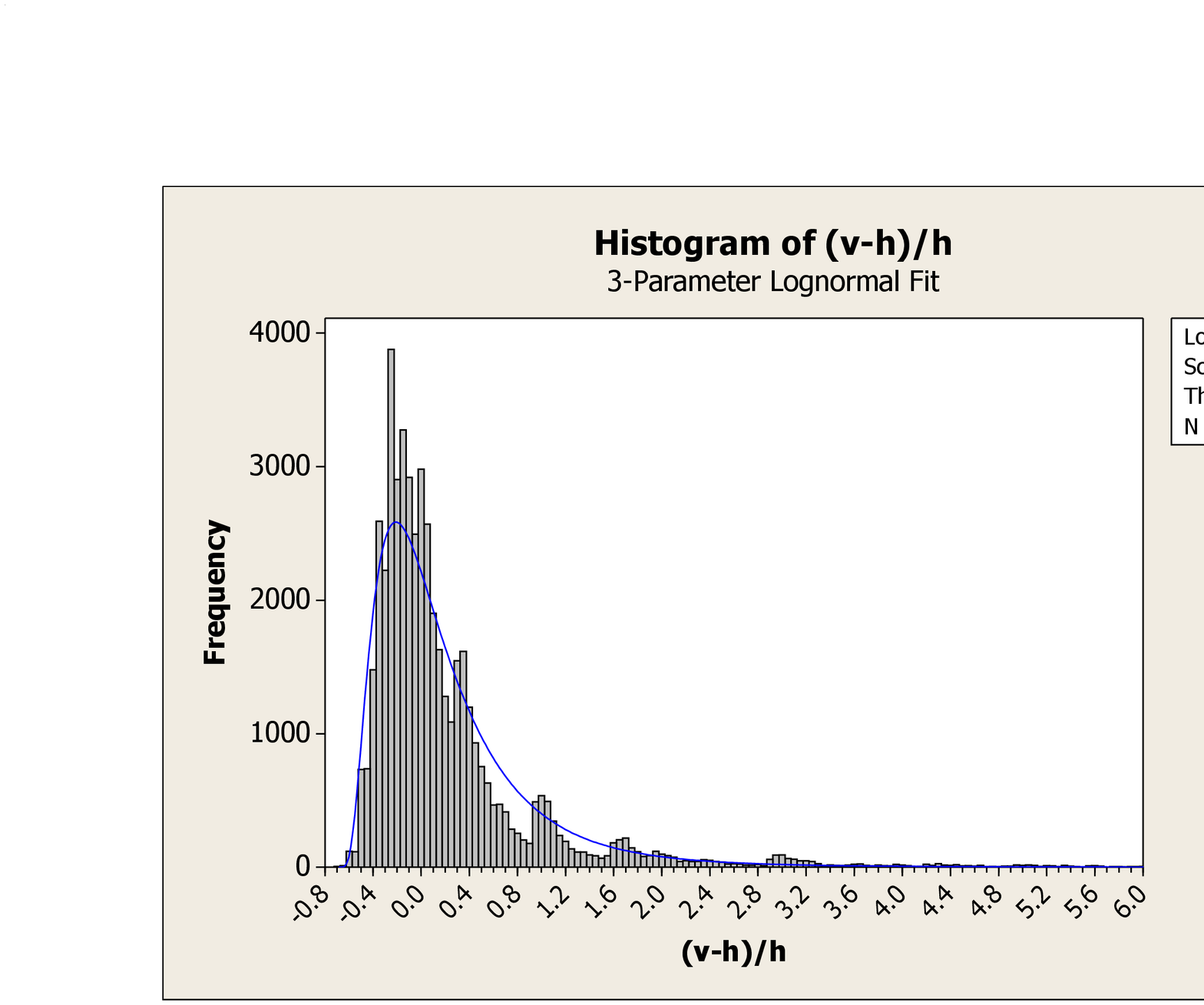} }

\centerline{{\bf Fig. 8.} Frequency histogram of $(v-h)/h$ with a
lognormal fit}
\centerline{$\min=-0.68,\ \max=14.77$ \quad (170 values omitted)}
\bigskip

\centerline{ \includegraphics*[scale=.5]{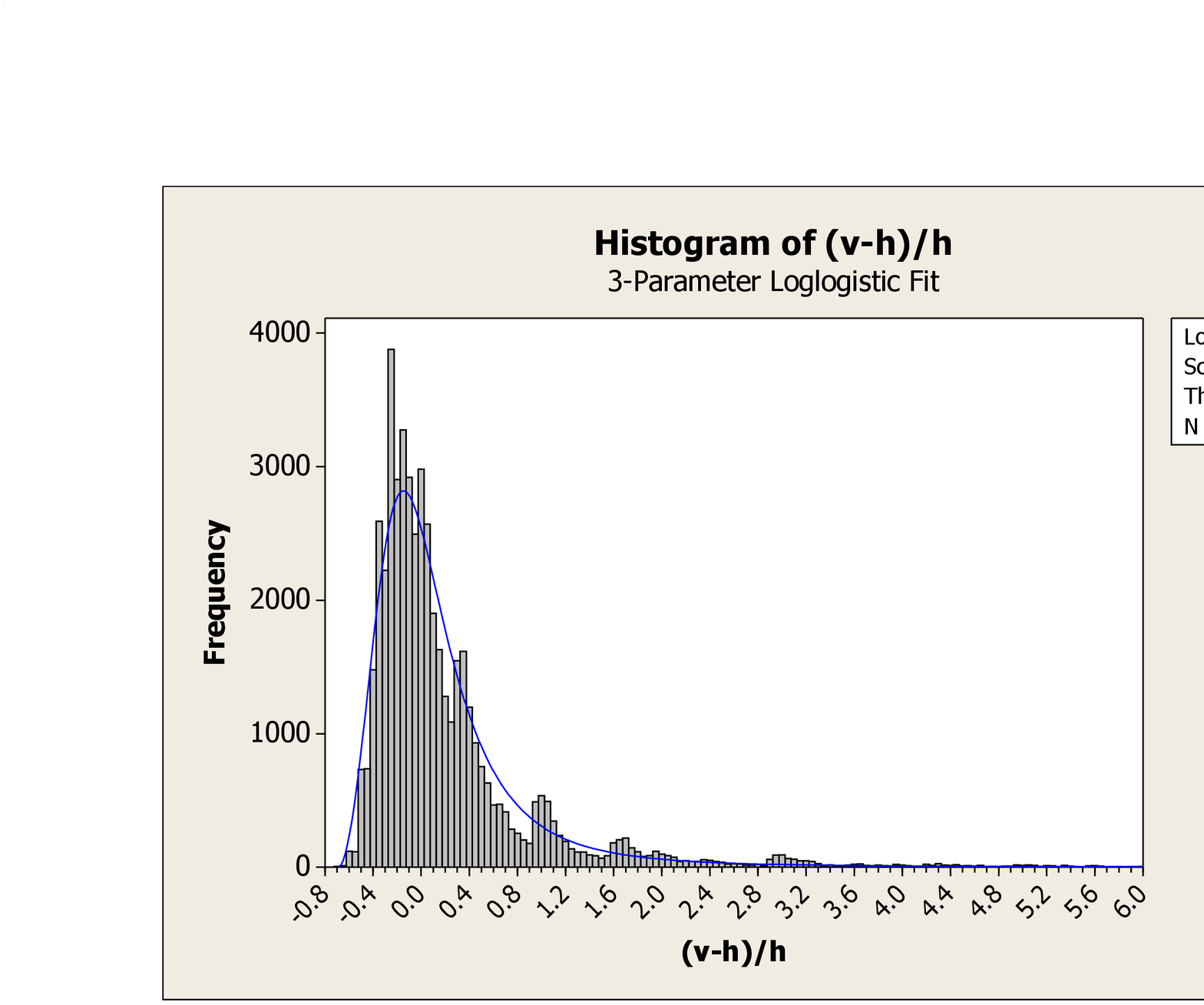} }

\centerline{{\bf Fig. 9.} Frequency histogram of $(v-h)/h$ with a
loglogistic fit}
\centerline{$\min=-0.68,\ \max=14.77$ \quad (170 values omitted)}
               \medskip

The histogram in Figures~6, 8, and~9 provides evidence that for most values of
$n$, $h(n)$ is a good approximation to $v(n)$. This leads to the main peak.
After comparing the histograms in Figures~6 and~7, it is plausible to
speculate that some of the secondary peaks of $(v(n)-h(n))$ to the
right of 0 correspond to large values of $\tau(n-1)$ that are quite
``popular''. It would be very
interesting to find (at least heuristically) a right model which
describes these secondary peaks (their height, frequency and so on).

Let $X$ be a random variable. We say that $X$ is {\em lognormally}
distributed if $\log X$ is  a normal distribution, and $X$ is  
{\em loglogistically} distributed if $\log X$ is a logistic distribution. The probability density functions of the lognormal distribution  is
$$ f(x;\mu,\sigma)= \frac{\exp(-(\log x-
\mu)^2/(2\sigma^2))}{\sqrt{2\pi}\sigma x}, $$
where $\mu$ and $\sigma^2$ are the mean and variance 
of $\log(X)$. The probability density
function of the loglogistic distribution  is
$$ f(x;\mu,\sigma)= \frac{\exp((\log x-\mu)/\sigma)}
{\sigma x(1+\exp((\log x-\mu)/\sigma))^2},
$$
where $\mu$ is the scale parameter and $\sigma$ is the shape parameter.

In Figures~8 and~9 we have provided the scaled histograms of $(v-h)/h$
with the lognormal fit and the loglogisitic fit respectively, as both of
them seem to be reasonable approximations.
Numerically,  the loglogistic fit seems to be better. However here is a heuristic argument (articulated by one of the referees) suggesting that the lognormal is more accurate. By the  Erd{\H o}s-Kac theorem~\cite[III.4.4, Theorem 8]{Ten3}, $\omega(s)$ 
is normally
distributed, and since
$\tau(s) = 2^{\omega(s)+O(1)}$ for most integers $s$, we conclude
that  $\log \tau(s)$ is also
normally distributed.   Given the connection between $v(n)$ and the 
divisor functions, it seems reasonable to believe that a lognormal distribution is more accurate.

As a curiosity, we also mention that
in the highly asymmetric histograms of Figures~6, 8 and~9 we still have
$v(n) < h(n)$ in 25057 out of 50000 cases. It would be interesting
to understand whether this is a coincidence, or whether there is some
regular effect behind this.

Our heuristic explanation for the difference between $V(N)$ and
$H(N)$ is as follows.
Overall, $G_n$ behaves as  a ``pseudorandom'' set,
but (as we observed in
Theorem~\ref{thm:Lower Bound}) there are some ``regular points'' on the
convex
closure arising from the divisors of $n-1$.  For a  typical integer $n$,
these points have little effect, but for exceptional values of $n$,
they make a substantial contribution to the value of
$v(n)$ which is sufficient to interfere with the ``pseudorandom''
behavior of $G_n$.
To see this, it is useful to recall that although for most integers we
have
$$\tau(n-1) = (\log n)^{\log 2 + o(1)} = h(n)^{\log 2 + o(1)} ,$$
see~\cite[Theorem~432]{HarWr},
on the average we have
$$
\sum_{n=2}^N \tau(n-1) \sim N \log N \sim \frac{3N}{8} H(N),
$$
see~\cite[Theorem~320]{HarWr}.
Therefore, the contribution of $2 \tau(n-1)$ from
the points on the curves $\alpha_1(n)$ and $\beta_1(n)$
(see Theorem~\ref{thm:Lower Bound}) is negligible compared
to $h(n)$ for almost
all $n$, but on average are of the same order as
$0.75H(N)$. Thus it is plausible to assert that the values of $H(N)$
reflect only the
``pseudorandom'' nature of $G_n$, whereas the contribution of
$2\tau(n-1)$ from the
curves $\alpha_1(n), \beta_1(n)$ reflect certain ``regular''
properties of the points of $G_n$.

\subsection{Weighted average contribution of divisors}

The lower bound of Theorem~\ref{thm:Lower Bound} takes into account only the
contribution from the divisors of $n-1$. It is plausible to assume that the
divisors of $jn-1$, with ``small'' $j\ge 2$, also give some regular
contribution to $v(n)$. This probably requires some completely new
arguments since the contribution from such divisors is certainly not additive.

Experimenting with some weighted averages involving $\tau(jn-1)$ for
``small'' values of $j$,
we have found that $g_1(n)$ and $g_2(n)$ where
\begin{eqnarray*}
g_1(n) & = & 2(\tau(n-1)-1) +  2 \sum_{j=2}^{\fl{\log n}}
j^{-3/2}\tau(jn-1) ,\\
g_2(n) &= &2(\tau(n-1)-1) + 2e\sum_{j=2}^{\fl{\log n}}
  e^{-j} \tau(jn-1) ,
\end{eqnarray*}
to be   ``reasonable'' numerical approximations to $v(n)$.

It is too early to make
any substantiated conjecture about the true contribution from the
divisors of $jn-1$ with $j\ge 2$.  Numerical experiments 
for a much broader range as well as some new ideas are needed.
Nevertheless, our calculation raises the
following question.

\begin{question}
\label{quest:Expansion}
Are there ``natural'' coefficients $c_j$, $j=2,3 \ldots$, and function $J(n)$, such that
if we define $g(n)$ to be 
$$
g(n) = 2\tau(n-1) +   \sum_{j=2}^{J(n)}
c_j \tau(jn-1),
$$
then we have
$$
V(N) \sim \frac{1}{N-1} \sum_{n=2}^N g(n)
$$
as $N\to \infty$?
\end{question}

Clearly, if $V(N) \sim C\log N$, then
the answer to Question~\ref{quest:Expansion} is positive, and one could then set $J(n)=2$ and determine the value of $c_2$ by ``reverse engineering''. However we are asking for coefficients $c_j$ and a function $J(n)$ that can be explained by some intrinsic reasons, provided such reasons exist!

\subsection{The difference $v(n)-2(\tau(n-1)-1)$}

Another computer experiment that we ran on our random set of 50000
integers was to check the values of the difference $v(n)-2(\tau(n-1)-1)$.
The histogram of our experiment is given in Figure~10.

\medskip
\centerline{ \includegraphics*[scale=.5]{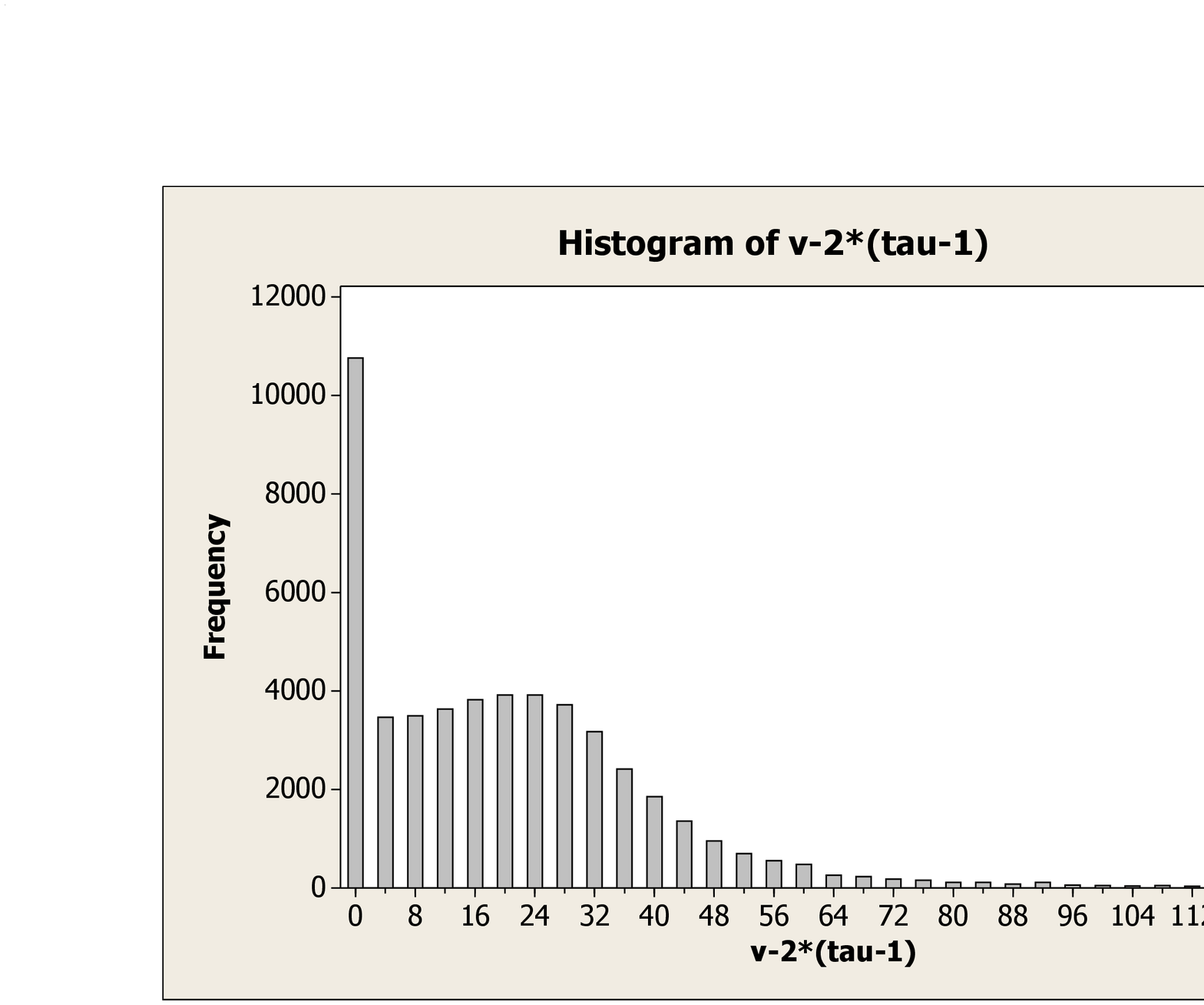} }

\centerline{{\bf Fig. 10.} Frequency histogram of $v(n)-2(\tau(n-1)-1)$}
\centerline{$\min=0,\ \max=484$ \quad (199 values omitted)}
\medskip

The  graph of Figure~10 suggests that the most ``popular'' value of
$v(n)-2(\tau(n-1)-1)$ is $0$. There is some obvious regularity in
the distribution of other values which would be interesting to explain.

The way we have derived the lower bound of Theorem~\ref{thm:Lower Bound} on the
frequency of the occurrence $v(n) = 2\( \tau(n-1)-1 \)$
from~\eqref{eq:m and T(n-1)} raises the following question:

\begin{question}
Is $T(n-1) =O(1)$
for all (or nearly all) integers $n$ with $v(n) = 2\( \tau(n-1)-1 \)$?
\end{question}

An affirmative answer to this question would then allow us to conclude that
$$
\# \left\{n \le x\ : \  v(n) = 2\( \tau(n-1)-1 \)\right\} \asymp
\frac{x}{\log x}.
$$

In our random set of 50000 integers we have 10764 integers satisfying
the equality
$v(n)=2(\tau(n-1)-1)$. For this set of 10764 integers we have
computed the value of
$t(n)$, where $t(n)= \fl{(T(n-1)+3)/4}$. We give this histogram in Figure~11.
We remark that for 7198 integers of this sample the value of
$t(n)$
is 1, and for
2413 integers of this sample the value of $t(n)$ is 2.
Thus for at least 9611 integers out of 10764 cases, we have
$\Gamma_n \cap \alpha_2(n) = \emptyset.$

\medskip
\centerline{ \includegraphics*[scale=.5]{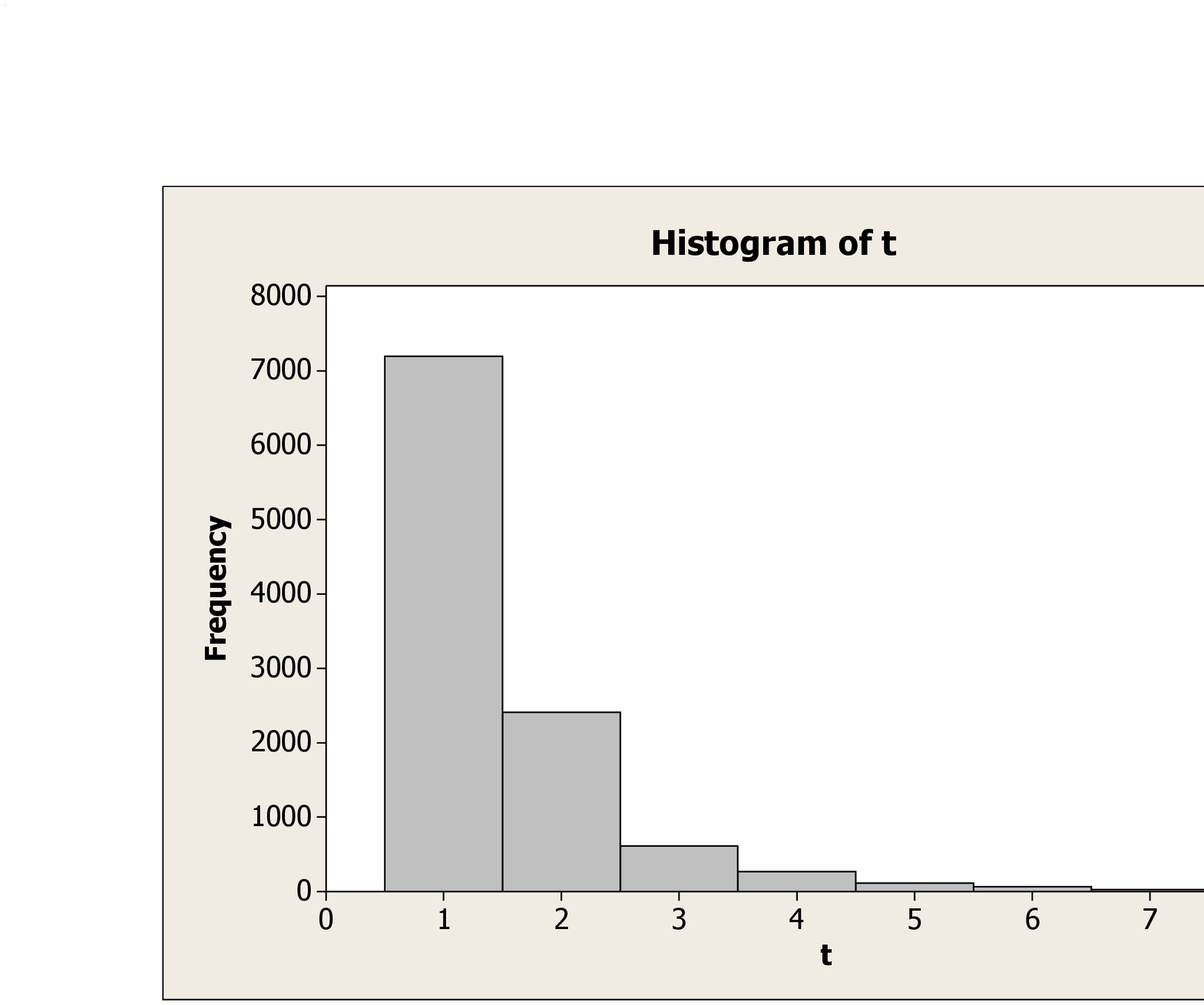} }

\centerline{{\bf Fig. 11.} Frequency histogram of $t(n)= \fl{(T(n-1)+3)/4}$}
\centerline{$\min=1,\ \max=26$ \quad (39 values omitted)}
\medskip

We have also found on examining the data that $v(n)-2(\tau(n-1)-1)$
is invariably a multiple of 4 and this suggests the following conjecture.

\begin{conjecture}
For almost all $n$,
$$
v(n)\equiv 2(\tau(n-1)-1) \pmod{4}.
$$
\end{conjecture}

We have a simple heuristic argument for this conjecture. We know that
$\tau(n-1)$ is odd if and only if $(n-1)$ is a square.
Thus the conjecture reduces to the statement that for almost all $n$,
$4\not|v(n)$. On invoking Propositions~\ref{Sym-Lem}
and~\ref{ver-cond} we have that $4|v(n)$ if and only if the vertex
$(a_s,b_s)$ lies on the line $x+y=n$. {\em Intuitively} this
seems to be a very rare occurrence (unfortunately at present we are
unable to put this key remark in a rigorous context); we
typically see that $a_s+b_s=n$ only when $n$ is the shifted square $m^2+1$.

\section{Other Curves}

Studying the point sets
$$
F_n(f) = \left\{  (a,b)\ : \ a,b\in \Z, f(a,b) \equiv 0
\pmod{n},\  0\leq a,b\leq n-1\right\},
$$
where $f(X,Y) \in \Z[X,Y]$,
is certainly a natural question, and this has been done in a number
of works, see~\cite{CZ,GrShZa,VZ,Zhe} and references therein.
In the case of prime modulus $p$,  one can use
the Bombieri~\cite{Bomb} bound of exponential sums along a curve
as a substitute of the bound of Kloosterman sums.
In particular, for a prime $n=p$, under some mild assumptions
on the polynomial $f$, one can easily obtain an analogue of
Theorem~\ref{thm:Upper Bound-1} for sets $F_p(f)$. However,
our other results are specific to the sets $G_n$ and cannot be extended
to other curves. It is worth remarking that for composite $n$, there are
some analogues of the Bombieri bound, see~\cite{SteShp}, but quite naturally,
they are much weaker than the bound of~\cite{Bomb}. So the
Kloosterman sums is one of very few examples where the
strength of the bound  remains almost unaffected
by the arithmetic structure of the modulus.

Our preliminary tests show that
the sets $F_n(f)$ and $F_p(f)$ have less ``infrastructure'' than
$G_n$ and behave more
like truly random sets.
For example, let $w_f(n)$ denote the
number of vertices of convex hull of $F_n(f)$. We now let
$$
h_f(n) =  \frac{8}{3}  \(\log \(\# F_n(f)\)  + \gamma  - \log 2 \).
$$
The histograms in Figures~12--14
show the relative difference $(w_f-h_f)/h_f$ for random quadratic and cubic
polynomials.
For the histogram of Figure~12 we chose a random value of $n$ in the
interval $[10000, 300000]$. Then based on the value of $n$ we
randomly chose the coefficients $a,b,c$ and took $f(x,y)$ to be the polynomial
$$ f(x,y)=y-ax^2-bx-c.$$
We did this for $10000$ values of $n$. For the histogram of Figure~13
we repeated
this same experiment with random quadratic polynomials for
$1000$ random primes in the interval $[7919,611953]$. For the
histogram of Figure~14 we  repeated  our first numerical experiment
(again for 10000 values of $n$), but this time with
random cubics
$$ f(x,y) = y-ax^3-bx^2-cx-d.$$

\centerline{ \includegraphics*[scale=.5]{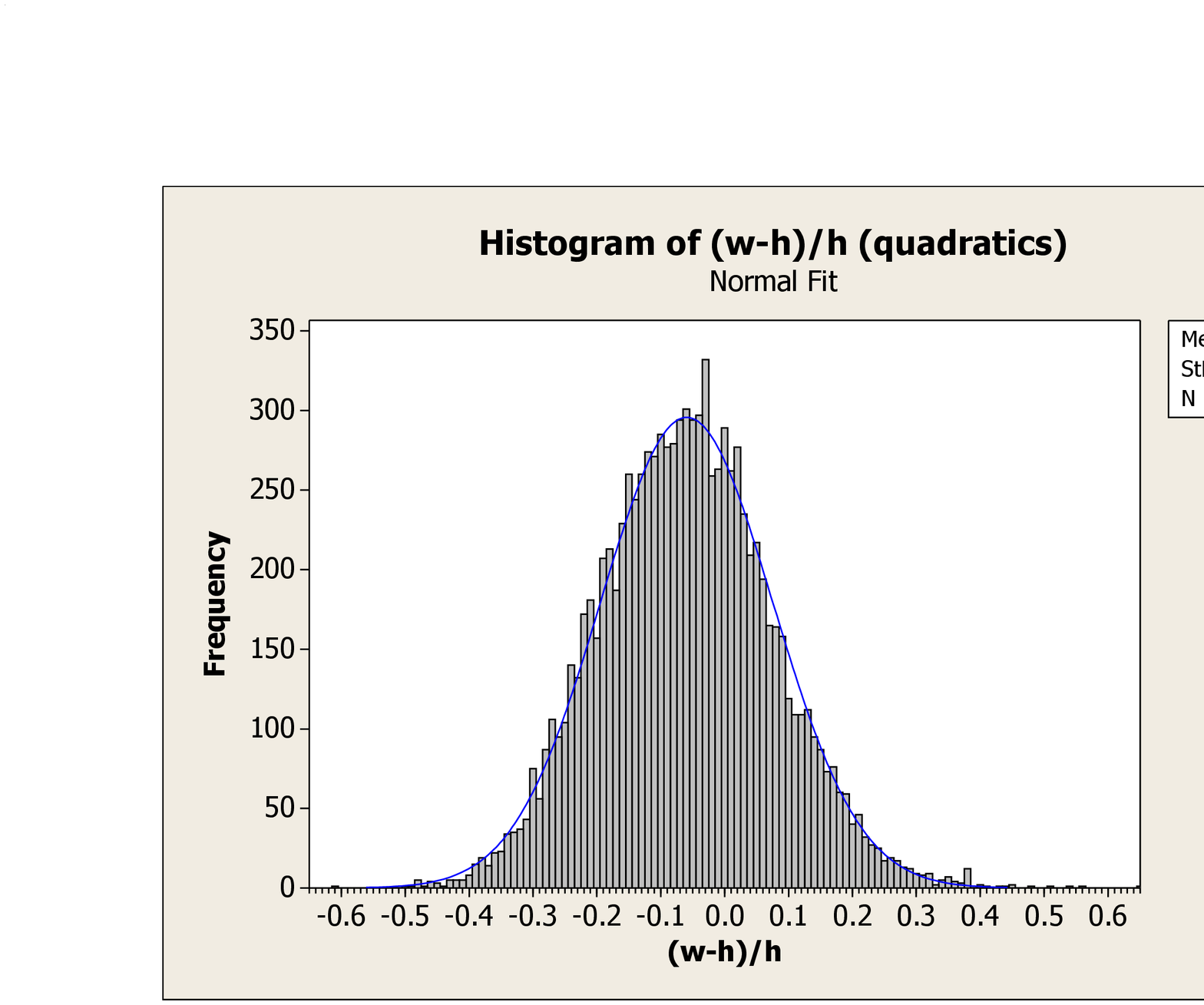} }

\centerline{{\bf Fig. 12.} Frequency histogram of $(w_f-h_f)/h_f$ for random
quadratics $f$ over random $n$}
\centerline{$\min=-0.607,\ \max=0.65$}
\bigskip

\centerline{ \includegraphics*[scale=.5]{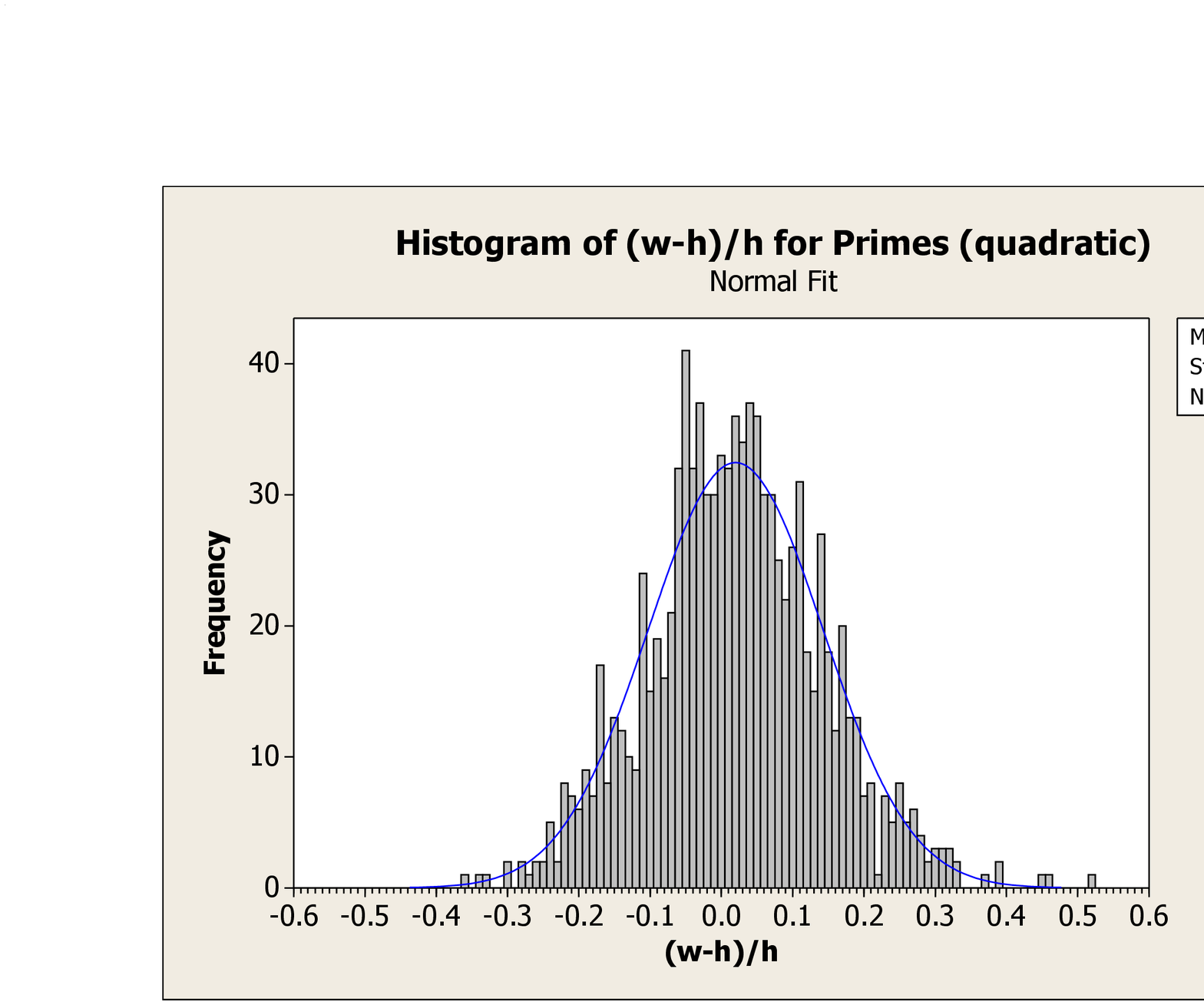} }

\centerline{{\bf Fig. 13.} Frequency histogram of $(w_f-h_f)/h_f$ for random
quadratics $f$ over random $p$.}
\centerline{$\min=-0.355,\ \max=0.518$}

\bigskip
\centerline{ \includegraphics*[scale=.5]{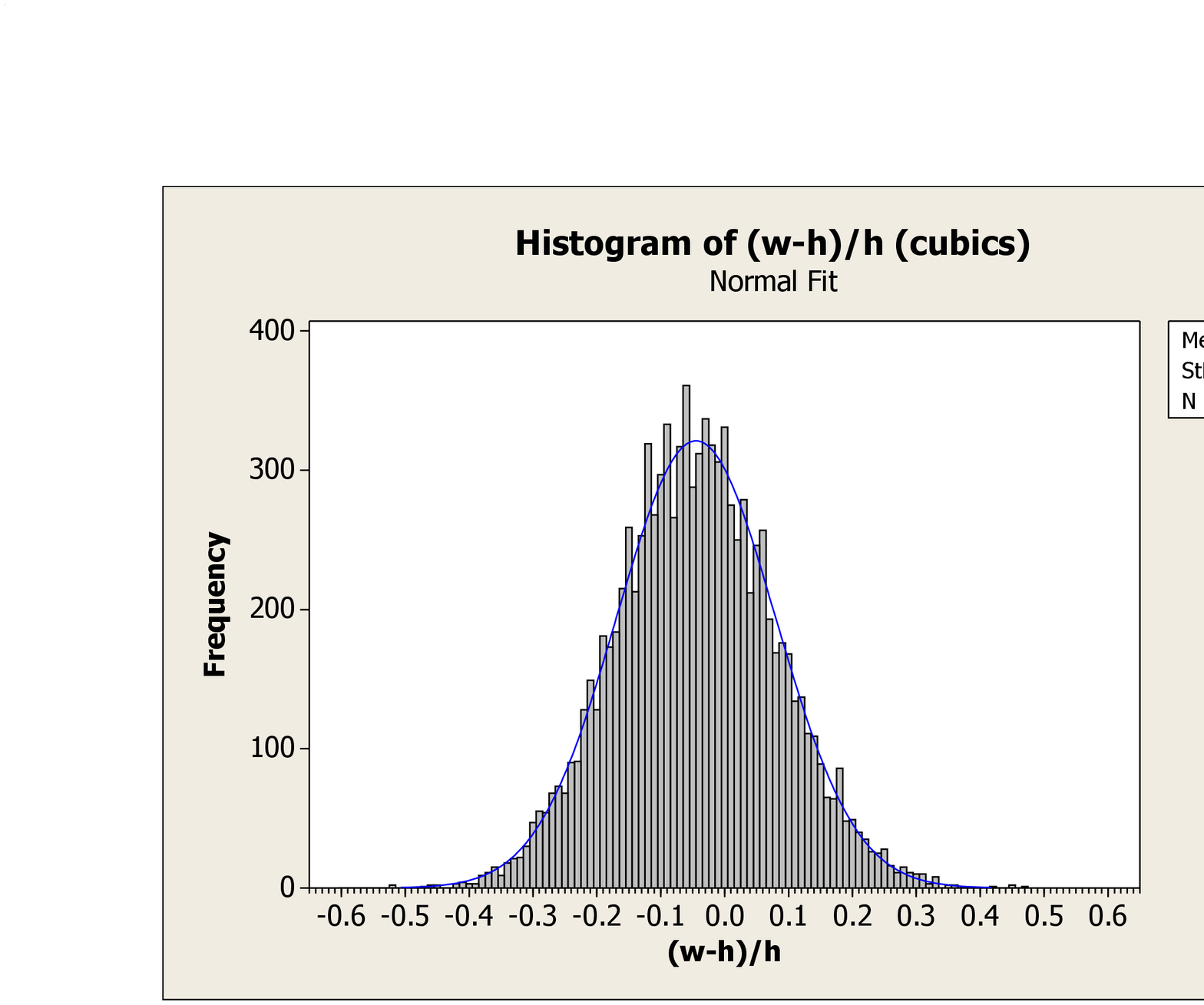} }

\centerline{{\bf Fig. 14.} Frequency histogram of $(w_f-h_f)/h_f$ for random
cubics $f$ over random $n$}
\centerline{$\min=-0.525,\ \max=0.473$}
               \medskip

The histograms of Figures~12--14 suggest that the quantities
$$\frac{w_f(n)-h_f(n)}{h_f(n)}  \qquad \text{and}\qquad
\frac{w_f(p)-h_f(p)}{h_f(p)}$$
are both normally distributed with mean 0, and so we make the
following ``Erd{\H o}s-Kac'' type   conjectures.

Let
$$
\Phi_\sigma(z) = \frac{1}{\sqrt{2\pi}\sigma}\int_{-\infty}^z
\exp\left (-\frac{t^2}{2\sigma^2} \right ) dt,
$$
denote the cumulative distribution function of a normal distribution with mean 0 and variance $\sigma^2$.

\begin{conjecture}
For each integer $n \ge 1$ we choose a sequence $\cF=(f_n)$ of polynomials
$f_n(x,y) \in \Z_n[x,y]$ of a fixed degree $d\ge 2$, chosen uniformly at random
over the residue ring $\Z_n$   and let
\begin{eqnarray*}
\sigma_\cF(N) &= & \sqrt{\frac{1}{N}\sum_{n \le N} \(w_{f_n}(n)/h_{f_n}(n)
-1 \)^2},\\
\rho_\cF(N) &= & \sqrt{  \frac{1}{\pi(N)}\sum_{p \le N}
\left(w_{f_p}(p)/h_{f_p}(p) -1 \right
)^2}.
\end{eqnarray*}
Then for any real $z$,
\begin{eqnarray*}
\frac{\# \left \{ n \le N \ : \ \(w_{f_n}(n)-h_{f_n}(n)\)/h_{f_n}(n)
\le z \right \}}{N\Phi_{\sigma_\cF(N)}(z)} & \to & 1 , \\
\frac{\# \left \{ p \le N \ : \ \(w_{f_p}(p)-h_{f_p}(p)\)/h_{f_n}(p)
\le z \right \}}{\pi(N)\Phi_{\rho_\cF(N)}(z)}
& \to  & 1,
\end{eqnarray*}
with probability $1$ (over the choice of $\cF=(f_n)$) as $N\to \infty$.
\end{conjecture}

\section*{Acknowledgements} We thank the following people:
\begin{itemize}
\item The referees for their careful reading of the article. The manuscript substantially benefited from their comments. In particular we are endebted to
the referee who suggested using the result of Saias~\cite{Saias} to
show that
$$\# \left\{n \le x\ : \  v(n) = 2\( \tau(n-1)-1 \)\right\} \gg 
\frac{x}{\log x}.$$
\item Kevin Ford for suggesting
the set $\cA(x)$ that arises in the proof of Theorem~\ref{thm:Strict Ineq}.
\item Daniel Sutantyo for computing $\psi(3/4)$.
\item  Anthony
Aidoo and Marsha Davis for assistance with the frequency histograms.
\end{itemize}

During the preparation of this paper,
I.~S.\ was supported in part by ARC grant DP0556431.

\end{document}